\documentclass{article}

\usepackage{epic,eepic,epsfig,amssymb,amsmath,amsthm}

\title{On the stabilization problem for nonholonomic distributions}

\author{L.~Rifford\footnote{Universit\'e de Nice-Sophia
    Antipolis, Labo.\ J.A.\ Dieudonn\'e, UMR 6621, Parc
    Valrose, 06108 Nice Cedex 02, France ({\tt
      rifford@math.unice.fr})}
\and
E.~Tr\'elat\footnote{Universit\'e d'Orl\'eans,
Math., Labo. MAPMO, UMR 6628,
Route de Chartres, BP 6759,
45067 Orl\'eans cedex 2, France ({\tt
  Emmanuel.Trelat@univ-orleans.fr})}}

\date{}

\newtheorem{theorem}{Theorem}

\newtheorem{proposition}[theorem]{Proposition}
\newtheorem{lemma}{Lemma}[section]
\theoremstyle{remark}\newtheorem{remark}{Remark}[section]
\theoremstyle{definition}

\theoremstyle{definition}

\def\R{\textrm{I\kern-0.21emR}}
\def\N{\textrm{I\kern-0.21emN}}

\newcommand{\SPAN}{{\rm Span}}

\newcommand{\co}{{\rm co}}

\newcommand{\ee}{\end{equation}}
\newcommand{\beq}{\begin{equation}}

\begin{document}

\maketitle

\begin{abstract}
Let $M$ be a smooth connected and complete manifold of dimension $n$,
and $\Delta$ be a smooth nonholonomic distribution of rank $m\leq n$
on $M$. We prove that, if there exists a smooth Riemannian metric on
$\Delta$ for which no nontrivial singular path is minimizing, then
there exists a smooth repulsive stabilizing section of $\Delta$ on
$M$. Moreover, in dimension three, the assumption  of
the absence of singular minimizing horizontal paths can be dropped in
the Martinet case. The proofs are based on the study, using specific
results of nonsmooth analysis, of an optimal control problem of Bolza
type, for which we prove that the corresponding value
function is semiconcave and is a viscosity solution of a
Hamilton-Jacobi equation, and establish fine properties of optimal
trajectories.
\end{abstract}

\section{Introduction}

Throughout this paper, $M$ denotes a smooth connected manifold of
dimension $n$.

\subsection{Stabilization of nonholonomic distributions}
Let $\Delta$ be a smooth distribution of rank $m\leq n$ on $M$, that
is, a rank $m$ subbundle of the tangent bundle $TM$ of $M$. This means
that, for every $x \in M$, there exist a neighborhood
$\mathcal{V}_{x}$ of $x$ in $M$, and a $m$-tuple $(f_{1}^x, \ldots,
f_{m}^x)$ of smooth vector fields on $\mathcal{V}_{x}$, linearly
independent on $\mathcal{V}_{x}$, such that 
$$
\Delta(y) =\SPAN
\left\{f_1^x(y),\ldots, f_m^x(y) \right\}, \quad \forall y\in \mathcal{V}_x.
$$
One says that the $m$-tuple of vector fields $(f_{1}^x,\ldots ,f_{m}^x)$
represents locally the distribution $\Delta$. The distribution
$\Delta$ is said to be \textit{nonholonomic} (also called totally nonholonomic \textit{e.g.} in \cite{asach}) if, for every $x \in M$,
there is a $m$-tuple $(f_{1}^x, \ldots, f_{m}^x)$
of smooth vector fields on $\mathcal{V}_{x}$ which represents locally the distribution and such that 
$$
\mbox{Lie} \left\{ f_1^x , \ldots, f_m^x \right\} (y) = T_yM, \quad \forall y \in \mathcal{V}_x, 
$$
that is, such that the Lie algebra spanned by $f_{1}^x, \ldots, f_{m}^x$, is equal to the whole tangent space $T_y M$, at every point $y\in \mathcal{V}_{x}$. This Lie algebra
property is often called \textit{H\"ormander's condition}.

An \textit{horizontal path} joining
$x_{0}$ to $x_{1}$ is an absolutely continuous curve  
$\gamma (\cdot) :[0,1]\rightarrow M$ such that $\gamma(0)=x_{0}$,
$\gamma(1)=x_{1}$, and such that
$\dot{\gamma}(t)\in\Delta(\gamma(t))$, for almost every
$t\in [0,1]$. According to the classical
Chow-Rashevsky Theorem (see \cite{bellaiche96,chow39,montgomery02,rashevsky38}),
since the distribution is nonholonomic on $M$,
any two points of $M$ can be joined by an horizontal path.

Let $\Delta$ be a nonholonomic distribution and $\bar{x} \in M$ be
fixed. We recall that, for a smooth vector field $X$ on $M$, the
dynamical system $\dot{x}=X(x)$ is said to be \textit{globally
  asymptotically stable} at the point $\bar{x}$, if the two following
properties are satisfied:
  \begin{itemize}
  \item[]\textit{Lyapunov stability:}  for every
neighborhood $\mathcal{V}$ of $ \bar{x}$, there exists a neighborhood
$\mathcal{W}$ of $ \bar{x}$ such that, for every $x\in \mathcal{W}$, the
solution of $\dot{x}(t)=X(x(t)), x(0)=x$, satisfies
$x(t) \in \mathcal{V}$, for every $t\geq 0$.
\item[]\textit{Attractivity:} for every $x\in M$, the solution of
      $\dot{x}(t)=X(x(t)),\ x(0)=x,$ tends to $\bar{x}$ as $t$ tends
      to $+\infty$.
\end{itemize}
The stabilization problem for nonholonomic
distributions consists in finding, if possible, a \textit{smooth stabilizing
  section} $X$ of $\Delta$, that is, a  smooth vector field $X$ on $M$
satisfying $X(x)\in \Delta(x)$ for every $x\in M$, such that the
dynamical system $\dot{x}=X(x)$ is globally asymptotically stable at
$\bar{x}$.

There exist two main obstructions for a distribution to admit a
stabilizing section. The first one is of global nature:
it is well-known that, if the manifold
$M$ admits such a dynamical system, then it possesses a smooth
Lyapunov function, \textit{i.e.}, a Morse function having only one
(possibly degenerate) critical point in $M$. In consequence, $M$ must
be homeomorphic to the Euclidean space $\R^n$ (we refer the reader to
\cite{riffordMCT03} for further details).
The second one is of local nature: due to \textit{Brockett's
  condition} (see \cite[Theorem 1, (iii)]{brockett}; see also
\cite{coron90,ryan94}),
the distribution $\Delta$ cannot admit a smooth stabilizing
section whenever $m<n$.

The absence of smooth stabilizing sections motivates to
define a new kind of stabilizing section. The first author has
recently introduced the notion of smooth repulsive stabilizing
feedback for control systems\footnote{If one represents locally the
  distribution $\Delta$ by a $m$-tuple of smooth vector
  fields $(f_1,\cdots,f_m)$, then the existence of a local stabilizing
  section for $\Delta$ is equivalent to the existence of a stabilizing
  feedback for the associated control system
  $\dot{x} = \sum_{i=1}^m u_i f_i(x)$. There is a large
  literature on alternative types of stabilizing feedbacks for control
  systems (see Section \ref{matin}).}
(see \cite{riffordMCT03,riffordTORINO,riffordSRS}), whose definition
can be easily translated in terms of stabilizing section.

Let $ \bar{x}\in M$ be fixed. Let $\mathcal{S}$ be a closed subset of $M$ and $X$ be a vector field on $M$. The
dynamical system $\dot{x}=X(x)$ is said to be \textit{smooth
repulsive globally asymptotically stable
at $ \bar{x}$ with respect to $\mathcal{S}$}
(denoted in short $\mbox{SRS}_{ \bar{x},\mathcal{S}}$) if the following
properties are satisfied:
\begin{itemize}
\item[(i)] The vector field $X$ is locally bounded on $M$ and smooth
  on $M\setminus \mathcal{S}$.
\item[(ii)] The dynamical system $\dot{x}=X(x)$ is globally
  asymptotically stable at $ \bar{x}$ in the sense of Carath\'eodory,
  namely, for every $x\in M$, there exists a solution of
\begin{eqnarray}\label{cara}
\dot{x}(t) =X(x(t)), \quad \mbox{for almost every } t \in [0,\infty),
\quad x(0)=x,
\end{eqnarray}
and, for every $x\in M$, every solution of (\ref{cara}) (called Carath\'eodory solution of $\dot{x}=X(x)$) on $[0,\infty)$
tends to $ \bar{x}$ as $t$ tends to $\infty$. Moreover, for every
neighborhood $\mathcal{V}$ of $ \bar{x}$, there exists a neighborhood
$\mathcal{W}$ of $ \bar{x}$ such that, for $x\in \mathcal{W}$, the
solutions of (\ref{cara}) satisfy
$x(t) \in \mathcal{V}$, for every $t\geq 0$.
\item[(iii)] For every $x\in M$, the solutions of (\ref{cara})
  satisfy
$
x(t) \notin \mathcal{S},
$
for every $t>0$.
\end{itemize}

In view of what happens whenever $\Delta =TM$, and having in mind the
above obstructions for the stabilization problem, a
natural question is to wonder if, given a smooth nonholonomic
distribution $\Delta$, there exists a section $X$ of $\Delta$ on $M$
and a closed nonempty subset $\mathcal{S}$ of $M$ such that $X$ is
$\mbox{SRS}_{ \bar{x},\mathcal{S}}$. In this paper, we provide a
positive answer in a large number of situations. To state our main
results, we need to endow the distribution $\Delta$ with a Riemannian
metric, thus encountering the framework of sub-Riemannian geometry,
and we require the concept of a singular path, recalled next.

%%%%%%%%%%%%%%%%%%%%%%%%%%%%%%%%%%%%%%%%%%%%%%%%%%%%%%%%%%%%%%%%%%%%%
\subsection{Sub-Riemannian geometry}

For $x_0\in M$, let $\Omega_\Delta(x_0)$ denote the set of horizontal
paths $\gamma(\cdot):[0,1]\rightarrow M$ such that $\gamma(0)=x_0$.
The set $\Omega_\Delta(x_0)$, endowed with the
$W^{1,1}$-topology, inherits of a Banach manifold
structure\footnote{It is
a straightforward adaptation of results of Bismut \cite{Bismut} (see also \cite{montgomery02}).}.
For $x_0,x_1\in M$, denote by $\Omega_\Delta(x_0,x_1)$ the set of
horizontal paths $\gamma(\cdot):[0,1]\rightarrow M$ such that
$\gamma(0)=x_0$ and $\gamma(1)=x_1$. Note that
$\Omega_\Delta(x_0,x_1)=\textrm{E}_{x_0}^{-1}(x_1)$, where
the end-point mapping $\textrm{E}_{x_0}:\Omega_\Delta(x_0)\rightarrow
M$ is the smooth mapping defined by $\textrm{E}_{x_0}(\gamma(\cdot))
:= \gamma(1)$.
A path $\gamma(\cdot)$ is said to be \textit{singular} if it
is horizontal and if it is a critical point of the end-point mapping
$\textrm{E}_{x_0}$.

The set $\Omega_\Delta(x_0,x_1)$ is a Banach submanifold of
$\Omega_\Delta(x_0)$ of codimension $n$ in a neighborhood of a
nonsingular path, but may fail to be a manifold in a neighborhood of
a singular path. It appears that singular paths play a crucial role
in the calculus of variations with nonholonomic constraints (see
\cite{CJT} for details and for properties of such curves).

Let $T^*M$ denote the cotangent bundle of $M$,
$\pi:T^*M\rightarrow M$ the canonical projection, and $\omega$ the
canonical symplectic form on $T^*M$. Let $\Delta^\perp$ denote
the annihilator of $\Delta$ in $T^*M$ minus its zero section. Define
$\overline{\omega}$ as the restriction of $\omega$ to
$\Delta^\perp$.
An absolutely continuous curve $\psi(\cdot):[0,1]\rightarrow
\Delta^\perp$ such that $\dot{\psi}(t) \in \ker
\overline{\omega}(\psi(t))$ for almost every $t\in [0,1]$, is
called an \textit{abnormal extremal} of $\Delta$.
It is well known that
a path $\gamma(\cdot):[0,1]\rightarrow M$ is singular if
and only if it is the projection of an abnormal extremal
$\psi(\cdot)$ of $\Delta$ (see \cite{Hsu} or \cite{CJT}).
The curve $\psi(\cdot)$ is said
to be an \textit{abnormal extremal lift} of $\gamma(\cdot)$.

\medskip

Let $g$ be a smooth Riemannian metric defined on the distribution
$\Delta$. The triple $(M,\Delta,g)$ is called a
\textit{sub-Riemannian} manifold.
The \textit{length} of a path $\gamma (\cdot) \in \Omega_{\Delta}(x_0)$ is
defined by
\begin{equation}\label{deflength}
\mbox{length}_{g}(\gamma (\cdot)) := \int_0^1
\sqrt{g_{\gamma(t)}(\dot{\gamma}(t),\dot{\gamma}(t))} dt.
\end{equation}
The \textit{sub-Riemannian distance} $d_{SR}(x_0,x_1)$ between two
points $x_0,x_1$ of $M$ is the infimum over the lengths
(for the metric $g$) of the horizontal paths joining $x_0$
and $x_1$. According to the Chow-Rashevsky Theorem (see
\cite{bellaiche96,chow39,montgomery02,rashevsky38}), since the distribution is
nonholonomic on $M$, the sub-Riemannian distance
is well-defined and continuous on $M\times M$. Moreover, if the
manifold $M$ is a complete metric space\footnote{Note that, since the
  distribution $\Delta$ is nonholonomic on $M$, the topology defined
  by the sub-Riemannian distance $d_{SR}$ coincides with the original topology of $M$ (see \cite{bellaiche96,montgomery02}).}
for the sub-Riemannian distance $d_{SR}$, then, since $M$ is
connected, for every pair $(x_{0},x_{1})$ of points of $M$
there exists an horizontal path $\gamma (\cdot)$ joining $x_{0}$ to $x_{1}$
such that
$$
d_{SR}(x_{0},x_{1}) = \mbox{length}_{g}(\gamma (\cdot)).
$$
Such an horizontal path is said to be \textit{minimizing}.

Define the Hamiltonian
$H:T^*M\rightarrow\R$ as follows. For every $x\in M$,
the restriction of $H$ to the
fiber $T^*_xM$ is given by the nonnegative quadratic form
\begin{equation}\label{monHam}
p \longmapsto \frac{1}{2}\max\left\{ \frac{p(v)^2}{g_x(v,v)}\ 
\vert\ v\in \Delta(x)\setminus\{0\} \right\} .
\end{equation}
Let $\overrightarrow{H}$ denote the Hamiltonian vector field on $T^*M$
associated to $H$, that is, $\iota_{\overrightarrow{H}}\omega=-dH$.
A \textit{normal extremal} is an integral curve of
$\overrightarrow{H}$ 
defined on $[0,1]$, i.e., a curve $\psi(\cdot):[0,1]\rightarrow T^*M$
such that $\dot{\psi}(t) = \overrightarrow{H}(\psi(t))$, for
$t\in[0,1]$. Note that the projection of a normal extremal is a
horizontal path.
The \textit{exponential mapping} $\textrm{exp}_{x_0}$ is defined on
$T^*_{x_0}M$ by $\textrm{exp}_{x_0}(p_0):=\pi(\psi(1))$, where
$\psi(\cdot)$ is the normal extremal so that $\psi(0)=(x_0,p_0)$ in
local coordinates.
Note that $H(\psi(t))$ is constant along a normal extremal
$\psi(\cdot)$, and that the length of the path $\pi(\psi(\cdot))$
is equal to $(2\, H(\psi(0)))^{1/2}$.

According to the Pontryagin maximum principle (see \cite{P}), a necessary
condition for a horizontal path to be minimizing is to be
the projection either of a normal extremal or of an abnormal
extremal. In particular, singular paths satisfy this
condition. However, a singular path may also be the projection of a
normal extremal.
A singular path is said to be \textit{strictly abnormal} if it is not
the projection of a normal extremal.

A point $x\in\textrm{exp}_{x_0}(T^*_{x_0}M)$ is said
\textit{conjugate} to $x_0$
if it is a critical value of the mapping $\textrm{exp}_{x_0}$. The
\textit{conjugate locus}, denoted by ${\mathcal{C}}(x_0)$,
is defined as the set of all points conjugate to $x_0$.
Note that Sard Theorem applied to the mapping $\textrm{exp}_{x_0}$
implies that the conjugate locus ${\mathcal{C}}(x_0)$ has Lebesgue
measure zero in $M$. 

\begin{remark}
It has been established in \cite{rt05} that the image of the
exponential mapping  $\textrm{exp}_{x_0}$ is dense in $M$, and is of
full Lebesgue measure for corank one distributions.
\end{remark}

\begin{remark}\label{remnormabnorm}
Let $x\in\textrm{exp}_{x_0}(T^*_{x_0}M)$, let $p_0\in T^*_{x_0}M$
such that $x= \textrm{exp}_{x_0}(p_0)$, and let $\psi(\cdot)$
denote the normal extremal so that $\psi(0)=(x_0,p_0)$ in local
coordinates. If $x$ is not conjugate to $x_0$,
then the path $x(\cdot):=\pi(\psi(\cdot))$ admits a unique normal extremal
lift. Indeed, if it had two distinct normal extremals lifts
$\psi_1(\cdot)$ and $\psi_2(\cdot)$, then the extremal
$\psi_1(\cdot)-\psi_2(\cdot)$ would be an abnormal
extremal lift of the path $x(\cdot)$. Hence, the path $x(\cdot)$ is
singular, and not strictly abnormal, and thus, in particular, the
point $x$ is conjugate to $x_0$. This is a contradiction.
\end{remark}

We also recall the notion of a cut point, required in this
article. Let $x_0\in M$;
a point $x\in M$ is not a \textit{cut point} with respect to
$ x_0$ if there exists a minimizing path joining $ x_0$ to $x$,
which is the strict restriction of a minimizing path starting
from $ x_0$.
In other words, a cut point is a point at which a minimizing
path ceases to be optimal.
The \textit{cut locus} of $ x_0$, denoted by
${\mathcal{L}}( x_0)$, is defined as the set of all cut points with
respect to $ x_0$. The following result is due to \cite{Sar}. We
provide in Section \ref{proofSar} a new (and selfcontained) proof of this
result, using techniques of nonsmooth analysis.

\begin{lemma}\label{lemconjcut}
Let $M$ be a smooth closed connected manifold of dimension $n$, and
$\Delta$ be a smooth nonholonomic distribution of rank $m\leq n$ on
$M$. Let $g$ be a metric on $\Delta$ for which no nontrivial singular
path is minimizing, and let $x_{0} \in M$. Then,
$$
\mathcal{C}_{min}(x_{0}) \subset \mathcal{L}(x_{0}),
$$
where $\mathcal{C}_{min}(x_{0})$ denotes the set of points $x\in M
\setminus \{x_{0}\}$ such that there exists
a critical point $p_0\in T_{x_{0}}^*M$ of the mapping
$\textrm{exp}_{x_0}$, and
such that the projection of the normal extremal $\psi (\cdot)$,
satisfying $\psi(0)=(x_0,p_0)$ in local coordinates,
is minimizing between $x_{0}$ and $x$.
\end{lemma}

In other words, under the assumptions of the lemma,
every (nonsingular) minimizing trajectory ceases to be
minimizing beyond its first conjugate point.

%%%%%%%%%%%%%%%%%%%%%%%%%%%%%%%%%%%%%%%%%%%%%%%%%%%%%%%%%%%%%%%%%%%%

\subsection{The main results}

\begin{theorem}
\label{THM1}
Let $M$ be a smooth connected manifold of dimension $n$,
and $\Delta$ be a smooth nonholonomic distribution of rank $m\leq n$
on $M$. Let $\bar{x}\in M$.
Assume that there exists a smooth Riemannian metric $g$ on
$\Delta$ for which $M$ is complete and no nontrivial singular path is minimizing. Then,
there exist a section $X$ of
$\Delta$ on $M$, and a closed nonempty subset $\mathcal{S}$ of $M$, of
Hausdorff dimension lower than or equal to $n-1$, such that $X$
is $\mbox{SRS}_{ \bar{x},\mathcal{S}}$.
\end{theorem}

\begin{remark}
If the manifold $M$, the distribution $\Delta$, and the metric $g$ are moreover real-analytic, then the
set $\mathcal{S}$ of the theorem can be chosen to be a subanalytic subset of $M \setminus \{\bar{x}\}$, of codimension greater than or equal to one (see \cite{Hardt,Hironaka}
for the definition of a subanalytic set). Note that, in this case,
since $\mathcal{S}$ is subanalytic (in $M \setminus \{\bar{x}\}$), it is a stratified (in the sense
of Whitney) submanifold of $M \setminus \{\bar{x}\}$.
\end{remark}

\begin{remark}
If $m=n$, then obviously there exists no singular path (it is the
Riemannian situation).
\end{remark}

\begin{remark}
The distribution $\Delta$ is called \textit{fat} (see
\cite{montgomery02}) at a point $x\in M$ if, for every vector field $X$
on $M$ such that $X(x)\in \Delta(x)\setminus\{0\}$, there holds
$$T_xM = \Delta(x)+\SPAN \{[X,f_i](x),\ 1\leq i\leq m\} ,  $$
where $(f_1,\ldots,f_m)$ is a $m$-tuple of vector fields representing
locally the distribution $\Delta$.
\\
With the same notations, it is called \textit{medium-fat} at $x$ (see
\cite{AS}) if there holds
$$T_xM = \Delta(x)+\SPAN\{[f_i,f_j](x),\ 1\leq i,j\leq m\}
+ \SPAN \{[X,[f_i,f_j]](x),\ 1\leq i,j\leq m\}.  $$

If $\Delta$ is fat at every point of $M$, then there exists no
nontrivial singular path (see \cite{montgomery02}). On the other part,
for a generic smooth Riemannian metric $g$ on $M$, every nontrivial
singular path must be strictly abnormal (see \cite{CJTnew}); it follows
from \cite[Theorem 3.8]{AS} that, if $\Delta$ is medium-fat at every
point of $M$, then, for generic metrics, there exists no nontrivial
minimizing singular path. Note that,
if $n\leq m(m-1)+1$, then the germ of a $m$-tuple of vector fields
$(f_1,\ldots,f_m)$ is generically (in $C^\infty$ Whitney topology)
medium-fat (see \cite{AS}).
\end{remark}

\begin{remark}
Let $m\geq 3$ be a positive integer, ${\cal G}_m$ be the set of pairs $(\Delta,g)$, 
where $\Delta$ is a rank $m$ distribution on $M$ and $g$
is a Riemannian metric on $\Delta$, endowed with the Whitney $C^\infty$
topology. There exists
an open dense subset $W_m$ of ${\cal G}_m$ such that every element
of $W_m$ does not admit nontrivial minimizing singular paths
(see \cite{CJTcras,CJT}).
This means that, for $m\geq 3$, generically, the main assumption of
Theorem \ref{THM1} is satisfied.
\end{remark}

In the following next result, we are able to remove,
in the compact and orientable three-dimensional case, the assumption on
the absence of singular minimizing paths.  Assume from now on that $M$
is a smooth closed manifold of dimension $3$ which is orientable  and
denote by $\Omega$ an orientation form on $M$. Any nonvanishing
one-form $\alpha$ generates a smooth rank-two distribution $\Delta$
defined by $\Delta := \ker\alpha $. Assume that $\Delta$ is
nonholonomic on $M$. There exists a unique smooth function $f$ on $M$
such that $\alpha \wedge d\alpha =f \Omega$ on $M$. Since $\Delta$ is
nonholonomic, the set $\{ f\neq 0 \}$ is open and dense in $M$. The
\textit{singular set} $\Sigma_{\Delta}$ of $\Delta$ is defined by
$$
\Sigma_{\Delta} := \left\{ x\in M \ \vert \ f(x) =0 \right\},
$$
Note that, if $M$ and $\alpha$ are analytic, then the singular set is an analytic subset of $M$.
The set $\Sigma_{\Delta}$ is said to be a \textit{Martinet
surface} if, for every $x\in \Sigma_{\Delta}$, $df(x) \neq 0$, so that the set $\Sigma_{\Delta}$ is 
a smooth orientable hypersurface on $M$. In the sequel, we will call a
\textit{Martinet distribution}, any nonholonomic distribution $\Delta$
associated with a nonvanishing one-form as above  such that
$\Sigma_{\Delta}$ is a Martinet surface.  In fact, it follows from the
generic classification of rank two
distributions on a three-dimensional manifold (see \cite{Zh}, see also
\cite{BTtoulouse}) that,
for every $x\in \Sigma_{\Delta}$, the distribution $\Delta$ is, in a
neighborhood of $x$, isomorphic to $\ker \alpha$, where the one-form
$\alpha$ is defined by $\alpha:=dx_3-x_2^2 dx_1$, in local coordinates
$(x_1,x_2,x_3)$. In this neighborhood, the Martinet surface
$\Sigma_\Delta$ coincides with the surface
$x_2=0$, and the singular paths are the integral curves of the vector
field $\frac{\partial}{\partial x_{1}}$ restricted
to $x_2=0$. This situation corresponds to the so-called
\textit{Martinet case}, and these singular paths are minimizing in the
context of sub-Riemannian geometry, for every smooth metric $g$ on
$\Delta$ (see \cite{ABCK,BTtoulouse,ls95}).

\begin{theorem}\label{THM2}
Let $M$ be a smooth connected orientable compact Riemannian manifold
of dimension three, and $\Delta$ be a Martinet distribution on
$M$. Let $\bar{x} \in M$. Then, there exist a section $X$ of
$\Delta$ on $M$, and a closed nonempty subset $\mathcal{S}$ of $M$, of
Hausdorff dimension lower than or equal to two, such that $X$
is $\mbox{SRS}_{ \bar{x},\mathcal{S}}$. 
\end{theorem}

\begin{remark}
The compactness assumption of the manifold $M$ can actually be dropped (see Remark \ref{remendproof2}). It is set to avoid technical difficulties in the proof.
\end{remark}

\subsection{Stabilization of nonholonomic control systems}\label{matin}
We begin this section with a remark on the local formulation of
Theorem \ref{THM1}. 
Let $U$ be an open neighborhood of $\bar{x}$ in $M$ such that
$\Delta_{\vert U}$ is spanned by a $m$-tuple $(f_1,\ldots,f_m)$ of
smooth vector fields on $U$, which are everywhere linearly independent
on $U$.
Every horizontal path $x(\cdot)\in\Omega(\bar{x})$, contained in $U$,
satisfies
\begin{equation}\label{grattecouilles}
\dot{q}(t)=\sum_{i=1}^mu_i(t)f_i(q(t))\ \textrm{for a.e.}\ t\in[0,1],
\end{equation}
where $u_i\in L^1([0,1],\R)$, for $i=1,\ldots,m$. The function
$u(\cdot)=(u_1(\cdot),\ldots,u_m(\cdot))$ is called the
\textit{control} associated to $x(\cdot)$, and the system
\ref{grattecouilles} is a \textit{control system}. Hence, Theorem
\ref{THM1}, translated in local coordinates, yields a stabilization
result for control systems of the form (\ref{grattecouilles}).

There are however slight differences between the geometric formulation
adopted in Theorem \ref{THM1}, and the corresponding result for
control systems. Indeed, when considering control systems of the form
(\ref{grattecouilles}), the vector fields $f_1,\ldots,f_m$ need not be
everywhere linearly independent. Moreover, a rank $m$ distribution
$\Delta$ on the manifold $M$ is not necessarily \textit{globally}
represented by a $m$-tuple of linearly independent vector fields (for
example, consider a rank two distribution on the two-dimensional
sphere of $\R^3$).

For these reasons, we derive hereafter a stabilization result, similar
to Theorem \ref{THM1}, valuable for control systems of the form
(\ref{grattecouilles}), and of independent interest.

\medskip

Consider on the manifold $M$ the control system
\begin{eqnarray}
\label{controlsystem}
\dot{x}(t) =\sum_{i=1}^m u_i(t) f_i(x(t)),
\end{eqnarray}
where $f_1,\ldots,f_m$ are smooth vector fields on $M$ (not
necessarily linearly independent), and the control 
$u=(u_1,\ldots,u_m)$ takes values in $\R^m$.

The system (\ref{controlsystem}) is said to be (totally) nonholonomic if the
$m$-tuple $(f_1,\cdots ,f_m)$ satisfies H\"ormander's condition
everywhere on $M$. According to the Chow-Rashevsky Theorem, any two
points of $M$ can be joined by a trajectory of (\ref{controlsystem}).

Let $\bar{x} \in M$ be fixed. The stabilization problem consists in
finding a \textit{feedback control function} $k=(k_1,\cdots ,k_m):
M\rightarrow \R^m$ such that the \textit{closed-loop system}
\begin{eqnarray}
\label{closedloop}
\dot{x} = \sum_{i=1}^l k_i(x) f_i(x)
\end{eqnarray}
is globally asymptotically stable at $\bar{x}$. It results from the
discussion above, and in particular from Brockett's condition, that
smooth or even continuous stabilizing feedbacks do not exist in
general. This fact has generated a
wide-ranging research with view to deriving adapted notions for
stabilization issues, such as discontinuous piecewise analytic feedbacks (see \cite{sussmann79}), discontinuous sampling feedbacks (see
\cite{clss97,rifford02}), continuous time varying control laws (see
\cite{coron92}), patchy feedbacks (see
\cite{anbr99}), almost globally asymptotically stabilizing feedbacks
(see \cite{riffordIHP} enjoying different
properties. The notion of smooth repulsive stabilizing feedback (see
\cite{riffordMCT03,riffordTORINO,riffordSRS}), whose definition is
recalled below, is under consideration in the present article.

Let $ \bar{x}\in M$ be fixed. Let $\mathcal{S}$ be a closed subset of $M$ and $k=(k_1,\cdots,k_m):M \rightarrow \R^m$ be a mapping on $M$. The feedback $k$ is said to be \textit{smooth repulsive globally asymptotically stable at $ \bar{x}$ with respect to $\mathcal{S}$} (denoted in short $\mbox{SRS}_{ \bar{x},\mathcal{S}}$) if the following
properties are satisfied:
\begin{itemize}
\item[(i)] The mapping $k$ is locally bounded on $M$ and smooth
  on $M\setminus \mathcal{S}$.
\item[(ii)] The dynamical system (\ref{closedloop}) is globally
  asymptotically stable at $ \bar{x}$ in the sense of Carath\'eodory.
\item[(iii)] For every $x\in M$, the Carath\'eodory solutions of (\ref{closedloop})
  satisfy
$
x(t) \notin \mathcal{S},
$
for every $t>0$.
\end{itemize}

We next associate to the control system (\ref{controlsystem}) an
optimal control problem.

For $x_0\in M$ and $T>0$\footnote{Note that, in what follows, the
  value of $T$ is not important. It can be assumed for instance that
  $T=1$.}, a control $u\in L^\infty([0,T],\R^m)$ is
said \textit{admissible} if the solution $x(\cdot)$ of
$(\ref{controlsystem})$ associated to $u$ and starting at $x_0$ is
well defined on $[0,T]$. On the set ${\cal U}_{x_0,T}$ of
admissible controls, and with the previous notations, define the {\em
  end-point mapping} by $E_{x_0,T}(u) := x(T)$.
It is classical that ${\cal U}_{x_0,T}$ is an open subset of
$L^\infty([0,T],\R^m)$ and that $E_{x_0,T}:{\cal U}_{x_0,T} \rightarrow
M$ is a smooth map.

A control $u\in {\cal U}_{x_0,T}$ is said to be {\em singular} if $u$
is a critical point of the end-point mapping $E_{x_0,T}$; in this case
the corresponding trajectory $x(\cdot)$ is said to be {\em singular}.

Let $x_0$ and $x_1$ be two points of $M$, and $T>0$.
Consider the optimal control problem of determining,
among all the trajectories of $(\ref{controlsystem})$ steering $x_0$
to $x_1$, a trajectory minimizing the {\em cost}
\begin{equation}\label{coutintro}
C_{U}(T,u)=\int_0^T u(t)^TU(x(t))u(t) dt,
\end{equation}
where $U$ takes values in the set $S_m^+$ of symmetric positive
definite $m\times m$ matrices.

\begin{theorem}\label{THM3}
Assume that there exists a smooth function $U:M\rightarrow S_m^+$ such
that no nontrivial singular trajectory of the control system
(\ref{controlsystem}) minimizes the cost (\ref{coutintro}) between its
extremities. Then, there exist a mapping $k:M\rightarrow \R^m$, and a  closed nonempty subset $\mathcal{S}$
of $M$, of Hausdorff dimension lower than or equal to $n-1$, such that $k$ is a $\mbox{SRS}_{ \bar{x},\mathcal{S}}$ feedback.
\end{theorem}

\begin{remark}
The same remarks as those following Theorem \ref{THM1} are
valuable. In particular, it is proved in \cite{CJTnew} that, for a
fixed smooth function $U:M\rightarrow S^+_m$, if $m\geq 3$, then
there exists an open and dense subset $O_{m}$ of the set of $m$-tuples
of smooth vector fields on $M$ so that the optimal control problem
(\ref{controlsystem})--(\ref{coutintro}) defined with an
$m$-tuple of $O_{m}$ does not admit nontrivial minimizing singular
trajectories.
\end{remark}

%%%%%%%%%%%%%%%%%%%%%%%%%%%%%%%%%%%%%%%%%%%%%%%%%%%%%%%%%%%%%%%%%%%%%%
%%%%%%%%%%%%%%%%%%%%%%%%%%%%%%%%%%%%%%%%%%%%%%%%%%%%%%%%%%%%%%%%%%%%%%

\section{Proof of the main results}
This section is organized as follows.
In Section \ref{secnonsmooth}, we recall some tools of nonsmooth
analysis that are required to prove our main results.
Section \ref{secbolza} is devoted to the proof of Theorem
\ref{THM1}. We first define a Bolza problem, equivalent to the
sub-Riemannian problem, for which we
derive some fine properties of the value function and of optimal
trajectories. In particular we prove that the value function is smooth
outside a singular set which is defined using a specific notion of
subdifferential. Theorem \ref{THM1} is then derived in Section
\ref{secproofthm1}. Theorem \ref{THM2} is proved in Section
\ref{secproofthm2}. The proof of Theorem \ref{THM3} is similar to the
one of Theorem \ref{THM1} and thus is skipped.

\subsection{Preliminaries: some tools of nonsmooth
  analysis}\label{secnonsmooth}
Let $M$ be a smooth manifold of dimension $n$.
\subsubsection{Viscosity subsolutions, supersolutions and solutions}
For an introduction to viscosity solutions of Hamilton-Jacobi
equations, we refer the reader to \cite{barles94,bcd97,CL,Lions}.
Assume that $F:T^*M \times \R \rightarrow \R$ is a continuous function on
$M$. A function $u:U\rightarrow \R$, continuous on the open set $U
\subset M$, is a \textit{viscosity subsolution}
(resp., \textit{supersolution}) on $U$ of
\begin{eqnarray}\label{HB}
F(x,du(x),u(x))=0,
\end{eqnarray}
if, for every $C^{1}$ function $\phi:U\rightarrow \R$ (resp., $\psi :U
\rightarrow \R$) satisfying $\phi \geq u$ (resp., $\psi \leq u$), and
every point $x_0 \in U$ satisfying $\phi(x_0) =u(x_0)$
(resp., $\psi(x_0)=u(x_0)$), there holds
$F(x_0,d\phi(x_0),u(x_0))\leq 0$
(resp., $F(x_0,d\psi(x_0),u(x_0))\leq 0$). A function is a
\textit{viscosity solution} of (\ref{HB}) if it is both a
viscosity subsolution and a viscosity supersolution of (\ref{HB}).

%%%%%%%%%%%%%%%%%%%%%%%%%%%%%%
\subsubsection{Generalized differentials}
Let $u:U \rightarrow \R$ be a continuous function on an open set $U
\subset M$. The \textit{viscosity subdifferential}
of $u$ at $x\in U$ is the subset of $T_x^*M$ defined by
$$
D^-u(x) := \left\{ d\psi (x) \ \vert \  \psi \in C^{1}(U) \mbox{ and }
  f-\psi \mbox{ attains a global minimum at } x\right\}.
$$
Similarly, the \textit{viscosity superdifferential}
of $u$ at $x$ is the subset of $T_x^*M$ defined by
$$
D^+u(x) := \left\{ d\phi (x) \ \vert \ \phi \in C^{1}(U) \mbox{ and }
  f-\phi \mbox{ attains a global maximum at } x\right\}.
$$
Notice that $u$ is a viscosity subsolution (resp., supersolution) of
(\ref{HB}) if and only if, for every $x\in U$ and every $\zeta \in
D^+u(x)$ (resp., $\zeta \in D^-u(x)$), one has $F(x,\zeta,u(x)) \leq
0$ (resp.,$F(x,\zeta,u(x)) \geq 0$).

The \textit{limiting subdifferential} of $u$ at $x\in U$ is the subset of $T_{x}^*M$ defined by
$$
\partial_{L} u(x) := \left\{ \lim_{k \rightarrow \infty} \zeta_{k} \ \vert \ \zeta_{k} 
\in D^-u(x_{k}), x_{k} \rightarrow x \right\}.
$$
By construction, the graph of the limiting subdifferential is closed
in $T^*M$. Moreover, the function $u$ is locally Lipschitzian on its
domain if and only if the limiting subdifferential of $u$ at any point
is nonempty and its graph is locally bounded (see \cite{clsw98,riffordmemoir}). 

Let $u:U \rightarrow \R$ be a locally Lipschitzian function.
The \textit{Clarke's generalized gradient} of
$u$ at the point $x\in U$ is the subset of $T_{x}^*M$ defined by
$$
\partial u(x) := \co \left(\partial_L u(x)\right),
$$
that is, the convex hull of the limiting differential of $u$ at $x$.
Notice that, for every $x\in U$,
$$
D^-u(x) \subset \partial_{L}u(x) \subset \partial u(x) \quad \mbox{and} \quad D^+u(x) \subset \partial u(x).
$$

\subsubsection{Locally semiconcave functions}
For an introduction to semiconcavity, we refer the reader to
\cite{cs04}. A function $u:U \rightarrow \R$, defined on the open set
$U \subset M$, is locally semiconcave on $U$, if for every $x \in U$,
there exist a neighborhood $U_x$ of $x$ and a smooth diffeomorphism
$\varphi_x: U_x \rightarrow \varphi_x(U_x) \subset \R^n$ such that
$f\circ \varphi_x^{-1}$ is locally semiconcave on the open subset $\tilde{U}_x =\varphi_x(U_x) \subset \R^n$. For the sake
of completeness, we recall that the function
$u: U \rightarrow \R$, defined on the open set $U \subset \R^n$, is
locally semiconcave on $U$, if  for every $\bar{x} \in U$
there exist  $C,\delta >0$ such that
\begin{eqnarray}
    \label{sc1}
\mu u(y)+ (1-\mu) u(x) -u(\mu x +(1-\mu) y) \leq \mu (1-\mu) C |x-y|^2,
\end{eqnarray}
for all $x,y \in \bar{x} + \delta B$ (where $B$ denotes the open unit
ball in $\R^n$) and every $\mu \in [0,1]$. This is equivalent
to say that the function $u$ can be written locally as
$$
u(x) = \left( u(x) - C|x|^2 \right) + \left( C|x|^2 \right), \quad \forall x \in \bar{x} + \delta B,
$$
that is, as the sum of a concave function and a smooth function. Note that every semiconcave function is locally Lipschitzian on its
domain, and thus, by Rademacher's Theorem, is differentiable almost
everywhere on its domain. The following result will be useful in the proof of our theorems.

\begin{lemma}
    \label{lem0}
    Let $u:U \rightarrow \R$ be a function defined on an open set
    $U\subset \R^n$. If, for every $\bar{x} \in U$, there exist a
    neighborhood $\mathcal{V}$ of $\bar{x}$ and a positive real number
    $\sigma$ such that, for every  $x\in \mathcal{V}$,
    there exists  $p_{x} \in \R^n$ such that  
    \begin{eqnarray}
        \label{sc2}
    u(y) \leq u(x) + \langle p_{x}, y-x\rangle + \sigma |y-x|^2,
    \end{eqnarray}
    for every  $y \in \mathcal{V}$, then the function $u$ is locally
    semiconcave on $U$.
    \end{lemma}
    
    \begin{proof}
        Without loss of generality, assume that $\mathcal{V}$ is an
        open ball $\mathcal{B}$.
     Let $x,y \in \mathcal{B}$ and $\mu \in [0,1]$. The point
     $\bar{x}:= \mu x + (1-\mu) y$ belongs to $\mathcal{B}$ by convexity. By assumption, there exists $\bar{p} \in \R^n$ such that  
     $$
     u(z) \leq u(\bar{x}) + \langle \bar{p}, z-\bar{x} \rangle + \sigma |z-\bar{x}|^2, \quad \forall z \in \mathcal{B}.
     $$
     Hence,
     \begin{eqnarray*}
     \mu u(y) + (1-\mu) u(x) & \leq & u(\bar{x}) + \mu \sigma |x-\bar{x}|^2 + (1-\mu) \sigma |y-\bar{x}|^2 \\
                             & \leq & u(\bar{x}) +  \left( \mu (1-\mu)^2 \sigma + (1-\mu) \mu^2 \sigma \right) |x-y|^2 \\
                             & \leq & u(\bar{x}) + 2 \mu (1-\mu) \sigma |x-y|^2,
        \end{eqnarray*}
        and the conclusion follows.
     \end{proof}

The converse result can be stated as follows.
 
 \begin{proposition}
 \label{propSC3}
Let $U$ be an open and convex subset of $\R^n$ and $u: U
\rightarrow \R$ be a function which is $C$-semiconcave on
$U$, that is, which satisfies
\begin{eqnarray}
    \label{SC33}
\mu u(y)+ (1-\mu) u(x) -u(\mu x +(1-\mu) y) \leq \mu (1-\mu) C |x-y|^2,
\end{eqnarray}
for every $x,y\in U$. Then, for every $x\in U$ and every $p\in D^+u(x)$, we
have
 \begin{eqnarray}
 \label{SC3}
 u(y) \leq u(x) + \langle p,y-x\rangle + \frac{C}{2} |y-x|^2, \quad \forall y \in \Omega,
 \end{eqnarray}
In particular, $D^+u(x) = \partial u(x),$ for every $x\in U$.
 \end{proposition}

\begin{remark}
\label{gena}
As a consequence (see \cite{cs04,riffordmemoir}), we obtain that, if a function $u:U \rightarrow \R$ is locally semiconcave
on an open set $U\subset M$, then, for every $x\in U$,
$$
\partial_L u(x) = \left\{ \lim_{k \rightarrow \infty} du (x_k) \ \vert \ x_{k} 
\in \mathcal{D}_u, x_{k} \rightarrow x \right\},
$$
 where $\mathcal{D}_u$ denotes the set of points of $U$ at which $u$ is differentiable.
 \end{remark}

The following result is useful to obtain several characterization of
the singular set of a given locally semiconcave function. We refer the
reader to \cite{cs04,riffordmemoir} for its proof.

\begin{proposition}
\label{propSC3*}
Let $U$ be an open subset of $M$ and $u: U
\rightarrow \R$ be a function which is locally semiconcave on
$U$. Then, for every $x\in U$, $u$ is differentiable at $x$
if and only if $\partial u(x)$ is a singleton.
\end{proposition}

The next result will happen to be useful (see \cite[Corollary
3.3.8]{cs04}).

\begin{proposition}
\label{propSCSC}
Let $u:U \rightarrow \R$ be a function defined on an open set
    $U\subset M$. If both functions $u$ and $-u$ are locally
    semiconcave on $U$, then $u$ is of class $C^{1,1}_{loc}$ on $U$.
\end{proposition}

%%%%%%%%%%%%%%%%%%%%%%%%%%%%%%

\subsubsection{Singular sets of semiconcave functions}

Let $u:U\rightarrow \R$ be a function which is locally semiconcave on
the open set $U\subset M$. We recall that since such a function is locally Lipschitzian on $U$, its limiting subdifferential is always nonempty on $U$. We define the \textit{singular set}
of $u$ as the subset of $U$
$$
\Sigma(u):= \left\{ x\in U \ \vert \ \partial_{L}u(x) \mbox{ is not a singleton} \right\}.
$$
Alberti, Ambrosio and Cannarsa proved in \cite{aac92} the following result.\footnote{In fact, this result has been strengthened later as follows. We can prove that the singular set of a locally semiconcave function is countably $n-1$-rectifiable, \textit{i.e.}, is contained in a countable union of locally Lipschitzian hypersurfaces of $M$ (see \cite{cs04,riffordmemoir}).}

\begin{theorem}
\label{THMAAC}
Let $U$ be an open subset of $M$. The singular set of a locally
semiconcave function $u:U \rightarrow \R$ is of Hausdorff
dimension lower than or equal to $n-1$.
\end{theorem}
 
The following lemma, proved in Appendix (Section
\ref{proofofLASTlem}), will be useful for the proof of Theorems
\ref{THM1} and \ref{THM2}.
 
\begin{lemma}
\label{LASTlem}
Let $u: U \rightarrow \R$ be a locally semiconcave function on an open
subset $U \subset M$ and  $\gamma :[a,b] \rightarrow U$ be a
locally Lipschitzian curve on the interval $[a,b]$. Then, for every
measurable map $p :[a,b] \rightarrow T^*M$ verifying
\begin{eqnarray*}
p(t) \in D^+  u(\gamma(t)), \quad \mbox{for a.e.} \ t \in [a,b],
\end{eqnarray*}
we have
\begin{eqnarray*}
\frac{d}{dt} \left( u(\gamma (t))\right) = p(t) \left( \dot{\gamma}(t)
\right),  \quad \mbox{for a.e.} \ t \in [a,b]. 
\end{eqnarray*}
\end{lemma}

%%%%%%%%%%%%%%%%%%%%%%%%%%%%%%%%%%%%%%%%%%%%%%%%%%%%%%%%%%%%%%%%%%%%

\subsection{Proof of Theorem \ref{THM1}}\label{secbolza}
From now on, assume that the assumptions of Theorem \ref{THM1}
hold. In particular, assume that there exists no nontrivial singular
minimizing path for the metric $g$.

\subsubsection{An equivalent optimal control problem}
Define the running cost $L_{g}$ by
$$
L_{g}(x,v) := g_{x}(v,v),
$$
for $x\in M$ and $v\in \Delta(x)$, and define
the functional $J_{g} : \Omega_{\Delta}(\bar{x}) \rightarrow \R^+$ by
$$
J_{g}(\gamma) := \int_{0}^1 L_g(\gamma(t),\dot{\gamma}(t)) dt.
$$
The Bolza optimization problem under consideration,
denoted by $(\mbox{BP})_{g,\Delta}$, consists in
minimizing the functional $J_g$, called \textit{energy}, over all
horizontal paths $\gamma$ joining $\bar{x}$ to $x\in M$.
Since $M$ is connected and complete, and since the running cost $L_g$
is coercive in every fiber, for every $x\in M$ there exists a
horizontal path $\gamma\in\Omega_\Delta(\bar{x},x)$, minimizing the energy
$J_g$. The value function associated to the Bolza problem
$(\mbox{BP})_{g,\Delta}$ is defined by
$$
V_{g,\Delta}(x) := \inf \left\{  J_g(\gamma) \ \vert \ \gamma \in
  \Omega_{\Delta}(\bar{x},x) \right\},
$$
for every $x\in M$.

\medskip

Note that the length of a horizontal path $\gamma$, defined by
(\ref{deflength}), does not depend on its parametrization. Hence, up
to reparametrizing, one can assume that the horizontal paths are
parametrized by arc-length, i.e., that
$g_{\gamma(t)}(\dot{\gamma}(t),\dot{\gamma}(t))=1$. In this case, the
length minimizing problem is equivalent to the minimal time
problem. Moreover, if all paths are defined on the same interval, then
length and energy minimization problems are equivalent, and
the value function $V_{g,\Delta}$ satisfies
\begin{equation}\label{eqdist}
V_{g,\Delta}(x) = d_{SR}(\bar{x},x)^2.
\end{equation}
In other terms, the sub-Riemannian problem of minimizing the length
between two points $\bar{x}$ and $x$, for the
sub-Riemannian manifold $(M,\Delta,g)$, is equivalent to the Bolza
problem $(\mbox{BP})_{g,\Delta}$.

\medskip

We next provide another equivalent formulation of this optimization
problem, in terms of optimal control theory, that will be useful in
the proofs of Theorems \ref{THM1} and \ref{THM2}.
Let $x\in M$, and let $\gamma$ be a minimizing horizontal path joining
$\bar{x}$ to $x$. Since $\gamma$ is necessarily not self-intersecting,
there exists a tubular neighborhood $\mathcal{V}$ of the path $\gamma$
in $M$, and there exist $m$ smooth vector fields $f_1,\ldots, f_m$ on
$\mathcal{V}$, such that
$$
\Delta (x) = \SPAN \left\{ f_i(x) \ \vert \ i=1,\ldots ,m
\right\},
$$
for every $x\in \mathcal{V}$.
Then, every horizontal path $x(\cdot)$, contained in $\mathcal{V}$, is
solution of the control system
$$ \dot{x}(t) = \sum_{i=1}^m u_i(t) f_i(x(t)),$$
where $u(\cdot) = (u_1(\cdot),\ldots,u_m(\cdot)) \in
L^{\infty}([0,1];\R^m)$ is called the \textit{control}.
Without loss of generality, we assume that the $m$-tuple of vector
fields $(f_1,\ldots, f_m)$ is orthonormal for the metric $g$.
In these conditions, the energy of the path $x(\cdot)$ is
$$
J_g(x(\cdot)) = \int_0^1 \sum_{i=1}^m u_i(t)^2 dt.
$$
Since the optimal control problem does not admit any nontrivial
singular minimizing path, it
follows from the Pontryagin maximum principle (see \cite{P})
that every minimizing
path $\gamma$ is the projection of a normal extremal
$\psi(\cdot)=(\gamma(\cdot),p(\cdot))$, associated with the control
 $u(\cdot) = (u_1(\cdot),\ldots,u_m(\cdot))$, where
\begin{equation}\label{contnormal}
u_i(t) = \langle p(t),f_i(\gamma(t))\rangle,\ \ i=1,\ldots,m.
\end{equation}

%%%%%%%%%%%%%%%%%%%%%%%%%%%%%%%%%%%

\subsubsection{Properties of the value function $V_{g,\Delta}$}

Consider the Hamiltonian function $H_{g,\Delta} : T^*M \rightarrow \R$
defined by
$$
H_{g,\Delta}(x,p) := \max_{v\in \Delta(x)} \left( p(v) - \frac{1}{2}
g_{x}(v,v) \right).
$$
Note that this Hamiltonian coincides with the Hamiltonian $H$ defined by (\ref{monHam}) (as can be 
seen in local coordinates).

\begin{proposition}
\label{PROP1}
If the distribution $\Delta$ is nonholonomic on $M$, then the value
function $V_{g,\Delta}:M \rightarrow \R$ is continuous on $M$ and is a
viscosity solution of the Hamilton-Jacobi equation
 \begin{equation}\label{HB1}
- \frac{1}{2}V_{g,\Delta} (x) + H_{g,\Delta} \left(x,\frac{1}{2}dV_{g,\Delta}(x)\right)=0, \quad \forall x \in M \setminus \{\bar{x}\}.
\end{equation}
\end{proposition}

Note that this proposition still holds if there exist some minimizing
singular paths.

\begin{proof}
    The continuity of $V_{g,\Delta}$ follows from the continuity of the sub-Riemannian distance, associated to the metric $g$, on 
    $M\times M$. Notice that, since the running cost $L_{g}$ is coercive in the fibers, and since $M$ is connected, 
    for every $x\in M \setminus \{\bar{x}\}$, there exists an horizontal 
    path $\gamma (\cdot) \in \Omega_{\Delta}(\bar{x},x)$ such that 
    $$
    V_{g,\Delta}(x) = J_{g}(\gamma(\cdot)) = \int_{0}^1 L_{g} (\gamma(s),\dot{\gamma}(s)) ds.
    $$
    Let us prove that $V_{g,\Delta}$ is a viscosity solution of (\ref{HB1}) on $M\setminus \{\bar{x}\}$. 
    Let $x\in M \setminus \{\bar{x}\}$, and let $\gamma:[0,1] \rightarrow M$ be an horizontal path joining $\bar{x}$ to $x$. 
    For $t\in (0,1)$, there exists $\tilde{\gamma} \in \Omega_{\Delta}(\bar{x},\gamma(t))$ such that 
\begin{eqnarray*}
V_{g,\Delta}(\gamma(t)) = \int_{0}^1 L_{g} (\tilde{\gamma}(s),\dot{\tilde{\gamma}}(s)) ds 
= \frac{1}{t} \int_{0}^t L_{g} \left(\tilde{\gamma}\left(\frac{s}{t}\right),\dot{\tilde{\gamma}}\left(\frac{s}{t}\right)\right) ds.
\end{eqnarray*}
Define $\gamma_{1} \in \Omega_{\Delta}(\bar{x},x)$ by
$$
\gamma_{1}(s) := \left\{
\begin{array}{ll}
    \tilde{\gamma} \left( \frac{s}{t}\right) & \mbox{ if } s \in [0,t], \\
    \gamma(s) & \mbox{ if } s \in [t,1].
    \end{array}
    \right.
$$    
Then, there holds
\begin{eqnarray*}
    V_{g,\Delta}(x) & \leq & \int_{0}^1 L_{g} (\gamma_{1}(s),\dot{\gamma}_{1}(s)) ds \\
    & \leq & \int_{0}^t L_{g} \left( \tilde{\gamma}\left(\frac{s}{t}\right), \frac{1}{t}\dot{\tilde{\gamma}}\left( \frac{s}{t}\right)\right) ds + \int_{t}^1 L_{g}(\gamma(s),\dot{\gamma}(s)) ds \\
    & \leq & \frac{1}{t} V_{g,\Delta}(\gamma(t)) + \int_{t}^1 L_{g} (\gamma(s),\dot{\gamma}(s)) ds.
\end{eqnarray*}
If $\phi:M \rightarrow \R$ is a $C^1$ function satisfying $\phi \geq V_{g,\Delta}$ and $\phi(x)= V_{g,\Delta}(x)$, 
then
\begin{eqnarray*}
    \phi(x) = V_{g,\Delta}(x) & \leq & \frac{1}{t} V_{g,\Delta}(\gamma(t)) + \int_{t}^1 L_{g} (\gamma(s),\dot{\gamma}(s)) ds \\
    & \leq & \frac{1}{t} \phi(\gamma(t)) + \int_{t}^1 L_{g} (\gamma(s),\dot{\gamma}(s)) ds.
 \end{eqnarray*}
Making $t$ tend to $1$, and considering all $C^1$ horizontal paths joining
$\bar{x}$ to $x$, we infer that, for every $v\in \Delta (x)$,
 $$
 d\phi(x) (v) \leq \phi(x)+ L_{g} (x,v).
 $$
On the other part, consider some path $\gamma \in
\Omega_{\Delta}(\bar{x},x)$ satisfying $ V_{g,\Delta}(x) =
J_{g}(\gamma)$. For every $t\in (0,1)$, up to a change of variable,
this path is necessary minimizing between $\bar{x}$ and $\gamma(t)$.
Therefore, for every $t\in (0,1)$,
$$
V_{g,\Delta}(x) = \frac{1}{t} V_{g,\Delta}(\gamma(t)) + \int_{t}^1
L_{g} (\gamma(s),\dot{\gamma}(s)) ds.
$$
If $\psi:U \rightarrow \R$ is a $C^1$ function satisfying $\psi \leq
V_{g,\Delta}$ and $\phi(x)= V_{g,\Delta}(x)$, then 
$$
\psi(x) \geq \frac{1}{t} \psi(\gamma(t)) +  \int_{t}^1 L_{g}
(\gamma(s),\dot{\gamma}(s)) ds.
$$
As previously, passing to the limit yields the existence of $v\in \Delta(x)$ such that
$$
d\psi(x)(v) \geq \psi(x) + L_{g}(x,v).
$$
The conclusion follows.
\end{proof}

\begin{remark}
\label{moule1}
Notice that, since $V_{g,\Delta}$ is a viscosity solution of
(\ref{HB1}) on $M\setminus \{\bar{x}\}$, there holds, for every horizontal path
$\gamma :[a,b] \rightarrow M\setminus \{\bar{x}\}$ (with $a<b$),
$$
V_{g,\Delta} (\gamma(b)) - V_{g,\Delta}(\gamma(a)) \leq \int_a^b V_{g,\Delta} (\gamma(s)) ds + \int_a^b L_g (\gamma(s),\dot{\gamma}(s)) ds.
$$
%We refer the reader to \cite{bcd97} for the proof of that result.
\end{remark}

\begin{remark}
\label{16juin1}
We also notice, that since $V_{g,\Delta}$  is a viscosity solution of
(\ref{HB1}) on $M\setminus \{\bar{x}\}$,  we have
\begin{equation}\label{HJBplus}
-\frac{1}{2}V_{g,\Delta}(x) + H_{g,\Delta} (x,\frac{1}{2}\zeta) =0,
\quad \forall x\in M\setminus \{\bar{x}\}, \quad \forall \zeta \in
\partial_L V_{g,\Delta}(x).
\end{equation}
\end{remark}

Finally we have the following result. 

\begin{proposition}\label{PROP2}
If the distribution $\Delta$ is nonholonomic on $M$,
then the value function $V_{g,\Delta}$ is continuous on $M$, and
locally semiconcave on $M \setminus \{\bar{x}\}$.
\end{proposition}

We just sketch the proof of Proposition \ref{PROP2}; we refer the
reader to \cite{cr04,riffordmemoir} for further details.

\begin{proof}
Recall that since $M$ is connected and complete, for every $x\in M \setminus \{\bar{x}\}$,
there exists an horizontal path $\gamma \in \Omega_{\Delta}(\bar{x},x)$
such that
$$
 V_{g,\Delta}(x) = J_{g}(\gamma) = \int_{0}^1 L_{g} (\gamma(s),\dot{\gamma}(s)) ds.
$$
By assumption, this minimizing path $\gamma$ is necessarily
nonsingular, and thus, it is the projection of a normal extremal. It
is well known (see \cite{agrachev98,trelat00}) that, for every $x\in M
\setminus \{\bar{x}\}$, there exists a neighborhood $\mathcal{V}$ of $x$
in $M \setminus \{\bar{x}\}$, such that the set of cotangent vectors
$p_0 \in T_{\bar{x}}^*M$ for which $\textrm{exp}_{\bar{x}}(p_0) \in
\mathcal{V}$ and the projection of the corresponding normal extremal
minimizes the length between $\bar{x}$ and $\textrm{exp}_{\bar{x}}(p_0)$,
is compact in $T_{\bar{x}}^*M$. On the other hand,  we know from
\cite[Proposition 4 p.~153]{rt05}, that, if $\zeta \in \partial_L
V_{g,\Delta}(x)$, then there exists a normal extremal $\psi:[0,1]
\rightarrow T^*M$ whose projection is minimizing between $\bar{x}$ and $x$
and such that $\psi(1)=(x,\frac{1}{2}\zeta)$. This proves that the
function $V_{g,\Delta}$ is locally Lipschitzian on $M\setminus
\{\bar{x}\}$.

Let $x \in M \setminus \{\bar{x}\}$, and let $\bar{\gamma}$ be a
minimizing horizontal path joining $\bar{x}$ to $x$. By assumption,
this path is nonsingular, and thus, it is not a critical point of the
end-point mapping $E_{\bar{x}}$.
Hence, there exists a submanifold $N$ of
$\Omega_\Delta(\bar{x})$, of dimension $n$, such that the mapping
$$
\begin{array}{rcl}
 \mathcal{E} :  N  & \longrightarrow & M \\
 \gamma(\cdot) & \longmapsto &  E_{\bar{x}}(\gamma(\cdot)) = 
  \gamma(1),
\end{array}
$$
is a local diffeomorphism, from a neighborhood of
$\bar{\gamma}(\cdot)$ in $N$, into a neighborhood
$\mathcal{W}$ of $x=\bar{\gamma}(1)$. We infer that, for every
$y \in \mathcal{W}$,
$$
V_{g,\Delta}(y) \leq J_g(\mathcal{E}^{-1}(y)).
$$
Since $J_g$ is smooth on the submanifold $N$, up to diffeomorphism, one can
put a parabola over the graph of $J_g$ on $N$, and thus, over the
graph of the function $V_{g,\Delta}$ at every $x \in M
\setminus \{\bar{x}\}$. The second-order term of this parabola
depends on the minimizing controls which are associated to the
points $x$. Using the compactness of the minimizers that we recalled
above,  we deduce that the function  $V_{g,\Delta}$ is locally semiconcave on $M
\setminus \{\bar{x}\}$.
\end{proof}

In the sequel, the singular set of $V_{g,\Delta}$, denoted
$\Sigma (V_{g,\Delta})$, is
$$
\Sigma(V_{g,\Delta}):= \left\{ x\in M \setminus \{\bar{x}\} \ \vert \ \partial_{L} V_{g,\Delta}(x) \mbox{ is not a singleton} \right\}.
$$
Recall that, since the function $V_{g,\Delta}$ is locally semiconcave
on $M\setminus \{\bar{x}\}$, its limiting subdifferential is nonempty
at any point of $M \setminus \{\bar{x}\}$ (see \cite{cs04}).

%%%%%%%%%%%%%%%%%%%%%%%%%%%%%%%%%

\subsubsection{Properties of optimal trajectories of
  $(\mbox{BP})_{g,\Delta}$}\label{proofSar}

%% Throughout this section, we assume that the Bolza problem
%% $(\mbox{BP})_{g,\Delta}$ does not admit any nontrivial singular
%% minimizing path. It follows from the Pontryagin maximum
%% principle (see \cite{P}), and from 

We stress that, due to the assumption of the absence of
singular minimizing path, every minimizing curve of the Bolza
problem $\mbox{(BP)}_{g,\Delta}$ is the projection of a normal
extremal, i.e., an integral curve of the Hamiltonian vector field
$\overrightarrow{H}$ defined by (\ref{monHam}),
associated with $H$.
In particular, every minimizing curve of $\mbox{(BP)}_{g,\Delta}$ is
smooth on $[0,1]$.

\begin{lemma}
    \label{lem1}
For every $x\in M \setminus \{\bar{x}\}$ and every $\zeta \in \partial_{L}
V_{g,\Delta}(x)$, there exists a unique normal extremal $\psi (\cdot) :[0,1]
\rightarrow T^*M$ whose projection $\gamma(\cdot) :[0,1] \rightarrow
M$ is minimizing between $\bar{x}$ and $x$, and such that
$\psi(1)=(x,\frac{1}{2}\zeta)$ in local coordinates. In addition,
$\psi(\cdot)$ is the unique (up to a multiplying scalar)
normal extremal lift of $\gamma (\cdot)$.
\end{lemma}

\begin{proof}
The first part of the statement is a consequence of \cite[Proposition
4 p.~153]{rt05}. Uniqueness follows from Cauchy-Lipschitz
Theorem. Uniqueness (up to a multiplying scalar) of the normal
extremal lift of $\gamma(\cdot)$ is a consequence of the assumption of
the absence of singular minimizing paths (see \cite[Remark 8
p.~149]{rt05}).
\end{proof}
 
\begin{lemma}
    \label{lem2}
Let $x\in M \setminus \{\bar{x}\}$ and $\gamma (\cdot):[0,1]\rightarrow M$ be a
minimizing curve of $(\mbox{BP})_{g,\Delta}$ such that $\gamma(1)=x$. Then, for
every $t\in (0,1)$, the curve $\tilde{\gamma}^t (\cdot):[0,1] \rightarrow M$
defined by $\tilde{\gamma}^t(s):=\gamma (st)$, for $s\in [0,1]$,
is the unique minimizing curve of $(\mbox{BP})_{g,\Delta}$ steering $\bar{x}$ to
$\gamma (t)$.
Moreover, $\tilde{\gamma}^t (\cdot)$ is the projection of the normal extremal $(\tilde{\gamma}^t (\cdot), \tilde{p}^{t}(\cdot) )$ in local coordinates, with $\tilde{p}^{t}(s)=tp(st)$ for every $s\in[0,1]$.
\end{lemma}

\begin{proof}
We argue by contradiction. If there is another horizontal curve
$\gamma_2 (\cdot):[0,1] \rightarrow M$ which minimizes the sub-Riemannian
distance between $\bar{x}$ and $\gamma(t)$, then there exists a nontrivial
minimizing path $x(\cdot)$, joining the points $\gamma(t)$ and
$\gamma(1)=x$, and having two distinct normal extremal lifts
$\psi_1(\cdot)$ and $\psi_2(\cdot)$. Then, the extremal
$\psi_1(\cdot)-\psi_2(\cdot)$ is an abnormal
extremal lift of the path $x(\cdot)$. Hence, the path $x(\cdot)$ is
singular and minimizing, and this contradicts our assumption.

We next prove that the adjoint vector associated to $\tilde{\gamma}^t (\cdot)$ is given by $\tilde{p}^{t}(s)=tp(st)$ for $s\in[0,1]$. In local coordinates, using the expression
(\ref{contnormal}) of normal controls, $\gamma(\cdot)$ is solution of the system
$$\dot{\gamma}(t)=\sum_{i=1}^n \langle p(t), f_i(\gamma(t))\rangle f_i(\gamma(t)),\quad\textrm{for a.e.}\ t\in[0,1].$$
Hence, $\tilde{\gamma}^t (\cdot)$ is solution of
$$\frac{d}{ds}\tilde{\gamma}^t(s)= t \sum_{i=1}^n \langle p(st), f_i(\tilde{\gamma}^t(t))\rangle f_i(\tilde{\gamma}^t(t)),\quad\textrm{for a.e.}\ s\in[0,1].$$
The conclusion follows.
\end{proof}

\begin{lemma}
    \label{lem3}
Any normal extremal $\psi (\cdot):[0,1]\rightarrow T^*M$ whose projection is
minimizing between $\bar{x}$ and $x\in M\setminus\{\bar{x}\}$ satisfies
$\zeta \in \partial_{L} V_{g,\Delta}(x)$, where
$\psi(1)=(x,\frac{1}{2}\zeta)$ in local coordinates.
\end{lemma}

\begin{proof}
Let $\psi (\cdot) :[0,1]\rightarrow T^*M$ be a normal extremal whose
projection $\gamma (\cdot)$ is minimizing between $\bar{x}$ and $x\in
M\setminus\{\bar{x}\}$. Since $V_{g,\Delta}$ is locally semiconcave on
$M\setminus \{\bar{x}\}$,
its limiting subdifferential is always nonempty on $M\setminus \{\bar{x}\}$. We infer from Lemmas \ref{lem1} and \ref{lem2} that, 
for every $t\in (0,1)$, there holds $\partial_{L}
V_{g,\Delta}(\gamma(t)) = \{\zeta(t)\}$, where
$\psi(t)=(x(t),\frac{1}{2t}\zeta(t))$ in local coordinates. Consider a sequence $(t_{k})$ of real numbers 
converging to $1$. Then, on the one part,  the sequence $(\psi(t_{k}))$  converges to $\psi(1)$, and on the other part, 
by construction of the limiting subdifferential, $\zeta=\zeta(1) \in
\partial_{L} V_{g,\Delta}(x)$.
\end{proof}

\begin{lemma}\label{lem4}
The following inclusion holds:
$$
\overline{\Sigma(V_{g,\Delta})} \setminus \Sigma(V_{g,\Delta})
\ \subset\ \mathcal{C}_{min}(\bar{x}) \cup \{\bar{x}\}.
$$
In particular, the set $\overline{\Sigma(V_{g,\Delta})}$ is of
Hausdorff dimension lower than or equal to $n-1$.
\end{lemma}
    
\begin{proof}
Let $x \in \overline{\Sigma(V_{g,\Delta})} \setminus
\Sigma(V_{g,\Delta})$ such that $x\neq \bar{x}$. By definition, the set
$\partial_{L} V_{g,\Delta}(x)$ is a singleton. 
Hence by Lemmas \ref{lem1} and  \ref{lem3}, there is a unique
minimizing path $\gamma (\cdot) \in \Omega_{\Delta}(\bar{x}x)$ and a
unique normal extremal $\psi (\cdot):[0,1] \rightarrow T^*M$ such that
$\gamma (\cdot)=\pi(\psi (\cdot))$; moreover, $\partial_{L} V_{g,\Delta}(x)
=\{\zeta\}$, where $\psi(1)=(x,\frac{1}{2}\zeta)$ in local
coordinates.
We argue by contradiction; if $x \notin \mathcal{C}_{min}(\bar{x})$,
then the exponential mapping
$\mbox{exp}_{\bar{x}}$ is not singular at $p_0$, where
$\psi(0)=(\bar{x},p_0)$ in local coordinates.
Furthermore, since
$x\in \overline{\Sigma(V_{g,\Delta})}$, there is
a sequence of points $(x_{k})$ in $\Sigma(V_{g,\Delta})$ which
converges to $x$. For every $k$, the set $\partial_{L}
V_{g,\Delta}(x_{k})$ 
admits at least two elements. Hence for every $k$, there are 
two distinct normal extremals $\psi_{k}^1(\cdot),\psi_{k}^2(\cdot):[0,1]
\rightarrow T^*M$ such
that their projections $\gamma_{k}^1(\cdot),\gamma_{k}^2(\cdot)$ are  minimizing
between $\bar{x}$ and $x_{k}$. Since the limiting
subdifferential of $V_{g,\Delta}$ is a singleton, the sequences
$(\psi_{k}^1(1)),(\psi_{k}^2(1))$ converge necessarily
to $\psi(1)$. Moreover, by regularity  of the Hamiltonian flow, the
sequences $(\psi_{k}^1(0)),(\psi_{k}^2(0))$ converge
necessarily to $\psi(0)$.
But the exponential mapping $\mbox{exp}_{\bar{x}}$ must be a local
diffeomorphism from a neighborhood of $p_0$
into a neighborhood of $\pi(\psi(1))$. This is a
contradiction. The second part of the lemma follows from the fact
that the singular set $\Sigma(V_{g,\Delta})$ is of Hausdorff dimension
lower than or equal to $ n-1$ (see Theorem \ref{THMAAC}), 
and of the fact that the set $\mathcal{C}_{min}(\bar{x})$ is contained in
$\mathcal{C}(\bar{x})$ which is of Hausdorff dimension
lower than or equal to $n-1$ (by \cite[Theorem 3.4.3]{federer69}).
\end{proof}

\begin{lemma}\label{lem4'}
The function $V_{g,\Delta}$ is of class $C^{1}$ on the open set $M
\setminus \left( \overline{\Sigma(V_{g,\Delta})} \cup \{\bar{x}\}
\right)$.
\end{lemma}
    
\begin{proof}
The set $\partial_L V_{g,\Delta} (x)$ is a singleton for every $x$ in
the set $M \setminus \left( \overline{\Sigma(V_{g,\Delta})}     \cup
  \{\bar{x}\} \right)$ which is open in $M$. From Remark \ref{gena} and the fact that $u$ is differentiable at some $x\in M\setminus \{\bar{x}\}$ if and only if $x\notin \Sigma(u)$, we infer that $V_{g,\Delta}$ is of class $C^1$ on the set $M
\setminus \left( \overline{\Sigma(V_{g,\Delta})} \cup \{\bar{x}\}
\right)$.
\end{proof}

\begin{lemma}\label{moule2}
Let $x\in M\setminus \{\bar{x}\}$ and $\bar{\gamma}(\cdot) :[0,1]
\rightarrow M$ be a minimizing curve of $(\mbox{BP})_{g,\Delta}$ such
that $\bar{\gamma}(1)=x$. Let $U_x$ be an open neighborhood of $x$ and
$\varphi_x :U_x \rightarrow \varphi_x(U_x) \subset \R^n$ be a smooth
diffeomorphism such that $V:=V_{g,\Delta} \circ \varphi_x^{-1}$ is a
locally semiconcave on the open subset $U:=\varphi_x (U_x) \subset
\R^n$. Let $t\in (0, 1)$ be such that $\bar{\gamma} (s) \in U_x $ for
every $s\in [t,1]$. Then there exist a neighborhood $\mathcal{W}_t$ of
$\bar{\gamma}(t)$ and   $\sigma(t) >0$ such that
\begin{eqnarray}
\label{moule3}
V(y) \geq V(\varphi_x (\bar{\gamma}(t))) + dV(\varphi_x (\bar{\gamma}(t)))(y-\varphi_x (\bar{\gamma}(t))) - \sigma(t) | y- \varphi_x (\bar{\gamma}(t))|^2, \quad \forall y \in \mathcal{W}_t.
\end{eqnarray}
\end{lemma}

\begin{proof}
Without loss of generality, we assume that $M=\R^n$, that
$\varphi_x$ is the identity, and that the closure of $U_x$ is a
compact subset of $M \setminus \{\bar{x}\}$. Set $x_s := \bar{\gamma}(s)$,
for every $s\in [t,1]$. Since $V=V_{g,\Delta}$ is locally semiconcave
on $M \setminus \{\bar{x}\}$, there exists $\sigma  \in \R$ such that
\begin{eqnarray}
\label{moule4b}
V(y) \leq V(x_s) + dV(x_s) (y-x_s) + \sigma |y-x_s|^2, \quad \forall y \in U, \quad \forall s \in [t,1].
\end{eqnarray}
The horizontal path $\tilde{\gamma} (\cdot) :[0,1-t] \rightarrow M$,
defined by 
$$
\tilde{\gamma}(s):= \gamma (1-s), \quad \forall s\in [0,1-t],
$$
is minimizing between $x$ and $\gamma (t)$. Hence, by assumption, it
is nonsingular, and thus, it is not a critical point of the end-point
mapping
$$
\begin{array}{rcl}
 E_t :   \Omega_{\Delta}(x)  & \longrightarrow & M \\
 \gamma(\cdot) & \longmapsto &  \gamma(1-t).
\end{array}
$$
Therefore, there exists a submanifold $N$ of $\Omega_{\Delta} (x)$ of
dimension $n$, such that the mapping
$$
\begin{array}{rcl}
 \mathcal{E}_t :  N  & \longrightarrow & M \\
 \gamma(\cdot) & \longmapsto &  E (\gamma(\cdot)),
\end{array}
$$
is a local diffeomorphism, from a neighborhood of $\tilde{\gamma} (\cdot)$ in $N$, into a neighborhood $\mathcal{W}_t$ of $\tilde{\gamma} (1-t) =x_t$. From Remark \ref{moule1}, we infer that, for every $y\in \mathcal{W}_t$, 
\begin{eqnarray}
\label{petasse1}
V(y) \geq V(x) - \int_t^1 V \left( \gamma_y (s)\right) ds - \int_t^1 L_g \left(  \gamma_y(s),\dot{\gamma}_y (s) \right) ds, 
\end{eqnarray}
where $\gamma_y (\cdot) : [t,1] \rightarrow M$ is defined by
$$
\gamma_y (s) :=  \mathcal{E}_t^{-1}(y)(1-s), \quad \forall s \in [t,1].
$$
By (\ref{moule4b}), we have
\begin{eqnarray}
\label{petasse2}
- \int_t^1 V \left( \gamma_y (s)\right) ds \geq - \int_t^1 \left( V(x_s)+ dV(x_s) (\gamma_y(s)-x_s) + \sigma |\gamma_y(s) -x_s|^2 \right) ds.
\end{eqnarray}
Moreover, since $\bar{\gamma}(\cdot)$ is minimizing between $\bar{x}$ and
$x$,
$$
V(x) = V(x_t) + \int_t^1 V(x_s) ds + \int_t^1 L_g(x_s,\dot{\bar{\gamma}}(s)) ds.
$$
Hence, from (\ref{moule4b}), (\ref{petasse1}) and (\ref{petasse2}), we
deduce that, for every $y\in \mathcal{W}_t$,
$$
V(y) \geq V(x_t) +  \Phi(y),
$$
where 
$$
\Phi (y) := \int_t^1 \left( L_g (x_s,\dot{\bar{\gamma}}(s)) - L_g (\gamma_y(s),\dot{\gamma}_y(s)) \right) ds - \int_t^1 \left( dV(x_s) (\gamma_y(s)-x_s) + \sigma  |\gamma_y(s) -x_s|^2 \right) ds.
$$
Since the mapping $\Phi_t : \mathcal{W} \rightarrow \R$ is smooth and
since $\Phi_t (x_t)=0$, a parabola can be put under the
graph of $V$ at $x_t$. This proves (\ref{moule3}).
\end{proof}

\begin{lemma}\label{moule4}
The following inclusion holds:
$$
\mathcal{C}_{min}(\bar{x}) \ \subset \ \overline{\Sigma(V_{g,\Delta})}.
$$
\end{lemma}

\begin{proof}
Let $x \in \mathcal{C}_{min}(\bar{x})$; note that, by definition of
$\mathcal{C}_{min}(\bar{x})$, one has $x\neq \bar{x}$ . We argue by
contradiction. If $x$ does not belong to
$\overline{\Sigma(V_{g,\Delta})}$, then $V_{g,\Delta}$ is $C^1$ in a
neighborhood of $x$. This means that there exist a neighborhood
$\mathcal{V}$ of $x$ and $t\in (0,1)$ such that for every $y\in
\mathcal{V}$, there is a minimizing curve of $(\mbox{BP})_{g,\Delta}$
such that $\bar{\gamma}(t)=y$. From the previous lemma and by
compactness of the minimizers, we deduce that the function
$-V_{g,\Delta}$ is locally  semiconcave on $\mathcal{V}$. Hence by
Proposition \ref{propSCSC}, $V_{g,\Delta}$ is $C^{1,1}_{loc}$ in
$\mathcal{V}$. Define
$$
\begin{array}{rcl}
 \Psi: \mathcal{V} & \longrightarrow & T_{\bar{x}}^*M \\
     y & \longmapsto     & \psi(0),
\end{array}
$$
where $\psi(\cdot): [0,1] \rightarrow TM$ is the normal extremal
satisfying $\psi(1)=(y,\frac{1}{2}dV_{g,\Delta}(y))$. This mapping is
locally Lipschitz on $\mathcal{V}$. Moreover by construction, $\Psi$
is an inverse of the exponential mapping. This proves that $p_0 :=
\Psi(x)$ is not a conjugate point. We obtain a contradiction.
\end{proof}

\begin{lemma}
\label{moule5}
Let $p_0 \in T_{\bar{x}}^*M$ be such that $H(\bar{x},p_0)\neq 0$. There exist a
neighborhood $\mathcal{W}$ of $p_0$ in $T_{\bar{x}}^*M$ and $\epsilon>0$
such that every normal extremal so that $\psi(0)= (\bar{x},p)$ (in local
coordinates) belongs to $\mathcal{W}$ is minimizing on the interval
$[0,\epsilon]$.
\end{lemma}

The proof of Lemma \ref{moule5} is postponed to the Appendix (Section
\ref{proofmoule5}).

\medskip

We are now ready to provide a proof for Lemma \ref{lemconjcut}.

\noindent \textit{Proof of Lemma \ref{lemconjcut}.}
For the sake of simplicity, we assume that $M=\R^n$, endowed with the
Euclidean metric.  We have to prove that $\mathcal{C}_{min}(\bar{x}) \
\subset \  \mathcal{L}(\bar{x})$. Let $y \in \mathcal{C}_{min}(\bar{x})$. We
argue by contradiction. Assume that $y$ does not belong to
$\mathcal{L}(\bar{x})$. This means that there exists a minimizing curve
$\gamma (\cdot)$ of $(\mbox{BP})_{g,\Delta}$ and $t_y \in (0,1)$ such
that $\gamma(t_y) =y$. Set $x:=\gamma(1)$, and let $\bar{t}$ be the
minimum of times $t\in (0,1)$ such that $\gamma (t) \notin
\overline{\Sigma(V_{g,\Delta})}$. We claim that $\bar{t} \in
(0,t_y]$. As a matter of fact, we know by Lemma \ref{moule4} that
$\gamma (t_y)=y \in \overline{\Sigma(V_{g,\Delta})}$. Moreover, from
Lemma \ref{moule5} and the absence of (nontrivial) singular minimizing
path, the mapping
$$
\begin{array}{rcl}
 \mathcal{W} & \longrightarrow & M \\
  p  & \longmapsto     &  \pi (\psi (\epsilon)),
\end{array}
$$
where $\psi(0)=(\bar{x},p)$, is injective. Hence from the Invariance of
Domain Theorem\footnote{The Invariance of Domain Theorem states that,
  for a topological manifold $N$, if $f:N \rightarrow N$
  is continuous and injective, then it is open. We refer the reader to
  the book \cite{bredon93} for a proof of that result.}, this mapping
is open. Which means that $V=V_{g,\Delta}  $ is necessarily of class $C^1$ on a
neighborhood of each $\gamma (s)$ with $s\in (0,\epsilon]$. We
conclude that $\bar{t} \in (0,t_y]$.\\
Set $\bar{x}:=\gamma(\bar{t})$ and $x_s:= \gamma (s)$ for every $s\in
[0,1]$. By local semiconcavity of $V$ (see Proposition
\ref{PROP2}), there exists a neighborhood $\mathcal{V}$ of $\bar{x}$ in
$M\setminus \{\bar{x}\}$ and $\sigma \in \R$ such that
\begin{eqnarray}
\label{fouffe}
V(z') \leq V(z) + \langle dV(z), z'- z \rangle + \sigma |z'- z|^2, \quad \forall z,z' \in \mathcal{V}.
\end{eqnarray}
Let $\bar{p} \in T_{\bar{x}}^*M$ such that
$\bar{x}=\mbox{exp}_{\bar{x}}(\bar{p})$. Since $V$ is of
class $C^1$ in a neighborhood of the curve $s \in (0,\bar{t})
\mapsto \gamma(s)$, there exists a neighborhood $\mathcal{W}'$ of
$\bar{p}$ in $ T_{\bar{x}}^*M$ such that
\begin{eqnarray*}
\forall p \in \mathcal{W}', \quad H(\bar{x},p) = H(\bar{x},\bar{p}) \Longrightarrow V(\mbox{exp}_{\bar{x}}(p)) =V(\bar{x}). 
\end{eqnarray*}
Thus, by (\ref{fouffe}), we have for every $p\in \mathcal{W}'$ satisfying $H(\bar{x},p) =H(\bar{x},\bar{p})$,
\begin{eqnarray}
\label{XXL}
\langle dV(\bar{x}), \mbox{exp}_{\bar{x}}(p) -\bar{x} \rangle \geq -\sigma \left| \mbox{exp}_{\bar{x}}(p) -\bar{x} \right|^2.
\end{eqnarray}
Furthermore, from Lemma \ref{lem3}, there exist a neighborhood
$\mathcal{V'}$ of $\bar{x}$ and $\sigma' >0$ such that
\begin{eqnarray}
\label{XXL2}
V(z) \geq V(\bar{x}) + \langle dV(\bar{x}, z-\bar{x} \rangle - \sigma' |z-\bar{x}|^2, \quad \forall z \in \mathcal{V}'.
\end{eqnarray}
Without loss of generality, assume that $\mathcal{V}'=\mathcal{V}$ and
$\sigma' =\sigma$. For every $p\in \mathcal{W}'$, set $x(p)
:=  \mbox{exp}_{\bar{x}}(p)$. By (\ref{XXL}) and (\ref{XXL2}), we
deduce that for every $p\in \mathcal{W}'$ satisfying $H(\bar{x},p)
=H(\bar{x},\bar{p})$ and for every $z \in \mathcal{V}$, we have
\begin{eqnarray*}
V(z) & \geq  & V(\bar{x}) + \langle dV(\bar{x}), z-  x(p)\rangle + \langle dV(\bar{x}), x(p)-\bar{x} \rangle -\sigma |z-\bar{x}|^2 \\
& \geq & V(\bar{x}) + \langle dV(\bar{x}), z- x(p) \rangle - \sigma | x(p)-\bar{x} |^2 -\sigma |z-\bar{x}|^2.
\end{eqnarray*}
In conclusion, by (\ref{fouffe}), we obtain that for every  $p\in \mathcal{W}'$ satisfying $H(\bar{x},p) =H(\bar{x},\bar{p})$ and every $z \in \mathcal{V}$, we have 
$$
\langle dV(\bar{x}),z-  x(p) \rangle -\sigma |  x(p) -\bar{x} |^2 - \sigma |z-\bar{x}|^2 \leq  \langle dV(x(p)), z-x(p) \rangle + \sigma |z-x(p)|^2.
$$
Hence, for every $p\in \mathcal{W}'$ satisfying $H(\bar{x},p)
=H(\bar{x},\bar{p})$ and every $z \in \mathcal{V}$,
\begin{eqnarray}
\label{fouffe2}
2 \sigma |z-x(p)|^2 + \langle dV(x(p))-dV(\bar{x}), z-x(p)\rangle + 2\sigma \langle x(p)-\bar{x}, z-x(p)\rangle + 2\sigma |x(p)-\bar{x}|^2 \geq 0.
\end{eqnarray}
Now, since $\bar{x}= \gamma(\bar{t})$ belongs to
$\overline{\Sigma(V)}$, we know by Lemma \ref{lem4}  that the
exponential mapping is singular at $\bar{p}$. Define the mapping $\Phi
: \R^n \times \R^n \rightarrow \R^n \times \R^n$ by
$$
\forall (z,p) \in \R^n \times \R^n, \quad \Phi (x,p) := \psi (x(1),p(1)),
$$
where $(x(\cdot),p(\cdot)) :[0,1] \rightarrow T^* M$ is the normal extremal satisfying $(x(0),p(0))=(x,p)$. Since $\Phi$  is a flow, its differential is always invertible. Hence there exist $P \in \R^n$ and  $Q \in \R^n \setminus \{0\}$ such that 
$$
D \Phi (\bar{x},\bar{p}) \cdot (0,P) = (0,Q).
$$
This means that there exist two continuous functions $\epsilon_1,
\epsilon_2 :\R \rightarrow \R^n$, and a mapping $\lambda \mapsto
p(\lambda) \in \mathcal{W}'$ such that, for every $\lambda$
sufficiently small, the following properties are satisfied:
\begin{itemize}
\item[(i)] $H(\bar{x},p(\lambda)) =H(\bar{x},\bar{p})$;
\item[(ii)] $V(x_{\lambda}) = V(\bar{x})$ where $x_{\lambda} := x(p(\lambda))$;
\item[(iii)] $x_{\lambda} =\bar{x}  +  \lambda^2 \epsilon_1 (\lambda)$;
\item[(iv)] $dV(x_{\lambda}) = dV(\bar{x}) + \lambda Q + \lambda^2 \epsilon_2 (\lambda)$.
\end{itemize}
From (\ref{fouffe2}), we deduce that, for every  $z \in
\mathcal{V}$,
$$
2 \sigma |z-x_{\lambda}|^2 + \lambda \langle Q , z-x_{\lambda} \rangle + \lambda^2 \langle \epsilon_2 (\lambda), z-x_{\lambda} \rangle + 2 \sigma \lambda^2 \langle \epsilon_1 (\lambda), z-x_{\lambda} \rangle + 2\sigma \lambda^4 |\epsilon_1(\lambda)|^2 \geq 0.
$$
We can apply this inequality for every $\alpha$ sufficiently small with  $z=x_{\lambda}-\alpha Q$. This yields
$$
2 \sigma \alpha^2 |Q|^2 -\lambda \alpha |Q|^2 + \lambda^2 \alpha \langle \epsilon_2(\lambda), -Q\rangle  + 2 \sigma \lambda^2 \alpha \langle \epsilon_1 (\lambda),-Q\rangle  + 2\sigma \lambda^4 |\epsilon_1 (\lambda)|^2 \geq 0,
$$
for every $\lambda, \alpha$ sufficiently small. Taking
$\alpha:=\lambda \sqrt{\lambda}$, we find a contradiction. $\Box$

\begin{lemma}\label{lem5}
There holds
$$
\overline{\Sigma(V_{g,\Delta})} = \mathcal{L} (\bar{x}) \cup \{\bar{x}\}.
$$
In particular, the cut locus is closed in $M \setminus \{\bar{x}\}$, and
is of Hausdorff dimension lower than or equal to $n-1$.
\end{lemma}

\begin{proof}
From Lemma \ref{lem1}, any point of $\Sigma(V_{g,\Delta})$ is joined
from $\bar{x}$ by several minimizing curves. Hence, from Lemma \ref{lem2},
any such point belongs to the cut locus $\mathcal{L}(\bar{x})$. From
Lemmas \ref{lem4} and \ref{lemconjcut}, we deduce that
$$
\overline{\Sigma(V_{g,\Delta})} \subset \mathcal{L}(\bar{x}) \cup
\{\bar{x}\}.
$$
If $x\in M\setminus \{\bar{x}\}$ does not belong to
$\overline{\Sigma(V_{g,\Delta})}$, then, from Lemma \ref{lem4'},
the function $V_{g,\Delta}$
is of class $C^{1}$ in a neighborhood $U$ of $x$.
Then, the continuous mapping
$$ \begin{array}{rcl}
 F : U & \longrightarrow & T^*M \\
     x & \longmapsto     & F(x) = -
     \overrightarrow{H}(x,\frac{1}{2} d V_{g,\Delta}(x))
\end{array}$$
is such that $F(x)=(\bar{x},p_0)$, with $\mathrm{exp}_{\bar{x}}(p_0)=x$.
This means that the exponential mapping $\mathrm{exp}_{\bar{x}}$ is a
homeomorphism from $F(U)$ into $U$, of inverse mapping
$\mathrm{exp}_{\bar{x}}$. In particular, it follows that
$x\notin \mathcal{L}(\bar{x})$. The fact that $\bar{x}$ belongs
to $\overline{\Sigma(V_{g,\Delta})}$ results from \cite[Theorem
1]{agrachev98}.
\end{proof}

\begin{remark}
Lemma \ref{lem5} asserts that the cut locus $\mathcal{L} (\bar{x})$ has
Hausdorff dimension lower than or equal to $n-1$. Recently, proving a
Lipschitz regularity property of the distance function to the cut locus,
Li and Nirenberg showed in \cite{ln05} that the $(n-1)$-dimensional
Hausdorff measure of the cut locus in the Riemannian framework is
finite. It would be interesting to study the regularity of the
distance function to the cut locus to obtain  such a result in the
sub-Riemannian case.
\end{remark}

\begin{lemma}\label{lem6}
The function $V_{g,\Delta}$ is of class $C^{\infty}$ on the open set
$M \setminus \overline{\Sigma(V_{g,\Delta})}$.
Moreover, if $\gamma:[0,1] \rightarrow M$ is a minimizing curve for
$\mbox{(BP)}_{g,\Delta}$, then $\gamma(t) \notin
\overline{\Sigma(V_{g,\Delta})}$, for every $t\in (0,1)$.
\end{lemma}
    
\begin{proof}
Let $\gamma:[0,1] \rightarrow M$ be a minimizing curve for
$\mbox{(BP)}_{g,\Delta}$. It follows from Lemmas \ref{lemconjcut} and
\ref{lem5} that $\gamma(t) \notin
\overline{\Sigma(V_{g,\Delta})}$, for every $t\in (0,1)$.

Let $x\in M \setminus \overline{\Sigma(V_{g,\Delta})}$, and let
$\gamma(\cdot)$ be a minimizing horizontal path joining $\bar{x}$ to
$x$. By assumption, $\gamma(\cdot)$ is necessarily nonsingular, and 
admits a unique normal extremal lift $\psi(\cdot):[0,1]\rightarrow
T^*M$. From Lemmas \ref{lemconjcut} and \ref{lem5}, the point $x$ is
not conjugate to $\bar{x}$, and hence,
the exponential mapping $\mathrm{exp}_{\bar{x}}$
is a (smooth) local diffeomorphism from a neighborhood of $p_0$ into
a neighborhood of $x$, where $\psi(0)=(\bar{x},p_0)$ in local
coordinates.
As recalled in the first section, the length of the path
$\gamma(\cdot)=\pi(\psi(\cdot))$ is equal to
$(2\, H(\psi(0)))^{1/2}$. Since $\gamma(\cdot)$ is minimizing, it is
also equal to $d_{SR}(\bar{x},x)$.
Then, using local coordinates, and from (\ref{eqdist}), there holds
$$ V_{g,\Delta}(x) = 2 H(\bar{x},(\mathrm{exp}_{\bar{x}})^{-1}(x)),$$
in a neighborhood of $x$
(see also \cite[Corollary 1 p.~157]{rt05}).
It follows that $V_{g,\Delta}$ is of class
$C^\infty$ at the point $x$.
\end{proof}

%%%%%%%%%%%%%%%%%%%%%%%%%%%%

\subsubsection{Conclusion: proof of Theorem
  \ref{THM1}}\label{secproofthm1}
Define $\mathcal{S} := \overline{\Sigma(V_{g,\Delta})}$. From Lemma
\ref{lem5}, there holds $\mathcal{S}= \mathcal{L}(\bar{x})\cup
\{\bar{x}\}$. We next define a section $X$ of $\Delta$, that is smooth outside
$\mathcal{S}$. To this aim, it is convenient to consider local
coordinates, and to express the problem in terms of optimal control. 
Let $x\in M\setminus\mathcal{S}$. In a neighborhood $U$ of $x$, one
has, in local coordinates,
$$\Delta = \SPAN \{ f_1,\ldots,f_m\},$$
where $(f_1,\ldots,f_m)$ is a $m$-tuple of smooth vector fields which is orthonormal for the metric $g$. We proceed as in \cite{rifford02}. 

Let $x \in M\setminus\bar{x}$ be fixed (of course, we set $X(x):=0$ if
$x=\bar{x}$), pick some $\zeta \in \partial_L V_{g,\Delta}(x)$.
Note that, since $V_{g,\Delta}$ is smooth outside the set
$\mathcal{S}$, one has $\zeta = dV_{g,\Delta}(x)$ whenever
$x\in M \setminus \mathcal{S}$.
Define the control $\tilde{u}(x)=(\tilde{u}_1(x),\cdots , \tilde{u}_m(x))$ by
\begin{equation}\label{couille2}
\tilde{u}_i(x) := \frac{1}{2} \zeta(f_i(x)), \quad \forall i=1, \cdots ,m.
\end{equation}
For $x\in M \setminus \mathcal{S}$, $\tilde{u}_i(x)=\frac{1}{2} \langle dV_{g,\Delta}(x),f_i(x)\rangle$ is the closed-loop form of the optimal control (\ref{contnormal}). For $x\in\mathcal{S}$, the expression of $\tilde{u}_i(x)$ depends on the choice of $\zeta \in \partial_L V_{g,\Delta}(x)$. Define 
\begin{equation}
\label{couilleX}
X(x):= -\sum_{i=1}^m \tilde{u}_i(x) f_i(x).
\end{equation}
Geometrically, $X(x)$ coincides with the projection of $-\frac{1}{2} \zeta$ onto $\Delta(x)$.
At the point $\bar{x}$, we set $X(\bar{x})=0$. This defines a vector field $X$ on $M$, which is smooth on $M\setminus\mathcal{S}$, but may be totally discontinuous on $\mathcal{S}$.

\medskip

We next prove that $X$ is $\mbox{SRS}_{\bar{x},\mathcal{S}}$.
Property (i) is obviously satisfied, but
properties (ii) and (iii) are not so direct to derive.

We first prove that every minimizing trajectory yields a Caratheodory solution of $\dot{x}=X(x)$. Let $x\in
M\setminus\bar{x}$ be fixed and $\gamma(\cdot) :[0,1] \rightarrow M$
be a minimizing curve of the Bolza problem $(\mbox{BP})_{g,\Delta}$
between $\bar{x}$ and $x$. It follows from the Pontryagin maximum
principle that $\gamma$ is the projection of a normal extremal expressed in local coordinates by $\psi(\cdot)=\left( \gamma(\cdot),p(\cdot)\right)$. Let $t\in (0,1)$;
from Lemma \ref{lem2}, the curve $\tilde{\gamma}^t
(\cdot):[0,1] \rightarrow M$
defined by $\tilde{\gamma}^t(s):=\gamma (st)$, for $s\in [0,1]$,
is the unique minimizing curve of $(\mbox{BP})_{g,\Delta}$ steering $\bar{x}$ to
$\gamma (t)$. Moreover, from Lemma \ref{lem2}, it is the projection of the normal extremal $\tilde{\psi}^t(\cdot) = \left( \tilde{\gamma}^t (\cdot), \tilde{p}^t (\cdot)  \right)$, where $\tilde{p}^t (\cdot) $ is defined by $\tilde{p}^t (s)=tp(st), $ for every $s\in [0,1]$. It then follows from Lemmas \ref{lem3} and \ref{lem6} that, along the curve $\gamma (\cdot)$,
$$
dV_{g,\Delta} (\gamma(t)) = 2tp(t), \quad \forall t \in (0,1).
$$

Therefore, $\gamma(\cdot)$ is solution of
$$
\dot{\gamma}(t) = \frac{1}{2t} \sum_{i=1}^m \Big( dV_{g,\Delta} (\gamma(t)) ( f_i(\gamma(t)) \Big) f_i(\gamma(t)), \quad \textrm{a.e.\ on}\ (0,1),
$$
in local coordinates along $\gamma(\cdot)$. This implies that the curve $x(\cdot):[0,\infty) \rightarrow M$ defined by
$$
x(t) := \gamma \left(e^{-t}\right), \quad \forall t \in (0,\infty),
$$
is a Carath\'eodory solution of $\dot{x}=X(x)$ such that $x(0)=\gamma(1)=x$. 

We next prove that any Carath\'eodory solution of $\dot{x}=X(x)$, $x(0)=x$, tends to $\bar{x}$ as $t$ tends to $+\infty$. Having in mind the minimizing properties (by construction) of the vector field $X$, it suffices actually to prove the following lemma.

\begin{lemma}\label{lemmepasrestersurS}
Let $x(\cdot)$ be any Carath\'eodory solution of $\dot{x}=X(x)$. Then, there does not exist a nontrivial interval $[a,b]$ such that $x(t)\in\mathcal{S}$ for every $t\in[a,b]$.
\end{lemma}

\begin{proof}
The proof goes by contradiction. Assume that there exist $\epsilon >0$ and a curve $x(\cdot):[0,\epsilon] \rightarrow M$ such that 
$$
\dot{x}(t) =X(x(t)), \quad \mbox{for almost every} \ t\in [0,\epsilon],
$$
and 
$$
x(t) \in \mathcal{S}, \quad \forall t \in [0,\epsilon].
$$
In local coordinates in a neighborhood of $x(0)=x$, one has
$$
\dot{x} (t) = X(x(t)) = - \frac{1}{2}\sum_{i=1}^m \zeta_t \left( f_i(x(t)\right) f_i(x(t)), \quad  \mbox{for almost every} \ t\in [0,\epsilon],
$$
where $\zeta_t \in \partial_L V_{g,\Delta}(x(t))$ for almost every $t\in [0,\epsilon]$. At this stage, we need to use Lemma \ref{LASTlem}, whose proof is provided in Appendix (Section \ref{proofofLASTlem}). According to this lemma, using
(\ref{couille2}) and the Hamilton-Jacobi equation (\ref{HJBplus}) satisfied by $V_{g,\Delta}$ (see Remark \ref{16juin1}), we deduce that, for almost every $t\in [0,\epsilon]$,
\begin{equation}
\label{midi0}
\frac{d}{dt} \left( V_{g,\Delta}(x(t))\right)  =  \zeta_t \left( \dot{x}(t)\right) 
 =  - \frac{1}{2}  \sum_{i=1}^m \Big( \zeta_t ( f_i(x(t))) \Big)^2 
 = - H_{g,\Delta}(x(t),\zeta_t)
 = - 2  V_{g,\Delta} (x(t)),
\end{equation}
since the Hamiltonian function $H_{g,\Delta}(x,p)$ is quadratic in $p$.
Therefore,
\begin{equation}
\label{midi}
V_{g,\Delta}(x(t)) = V_{g,\Delta}(x) e^{-2t}, \quad \forall t\in [0,\epsilon].
\end{equation}
Let $\gamma (\cdot) \rightarrow M$ be a minimizing curve of the Bolza problem $(\mbox{BP})_{g,\Delta}$ between $\bar{x}$ and $x(\epsilon)$. Define the horizontal path $\tilde{\gamma}(\cdot) :[0,1] \rightarrow M$ by
$$
\tilde{\gamma}(t) = \left\{ \begin{array}{rll} 
x(-\ln t) & \mbox{ if } & e^{-\epsilon} \leq t \leq 1  \\
\gamma \left( e^{\epsilon} t\right) & \mbox{ if } & 0 \leq t \leq e^{-\epsilon}.
\end{array}
\right.
$$
The cost of $\tilde{\gamma}(\cdot)$ is
\begin{eqnarray*}
J_{g}(\tilde{\gamma}(\cdot)) &= & \int_0^{e^{-\epsilon}} L_g \left( \tilde{\gamma}(t),\dot{\tilde{\gamma}}(t)\right) dt + \int_{e^{-\epsilon}}^1 L_g \left( \tilde{\gamma}(t),\dot{\tilde{\gamma}}(t)\right) dt \\
& = & \int_0^{e^{-\epsilon}} L_g\left(\gamma \left(e^{\epsilon}t\right),e^{\epsilon} \dot{\gamma} \left(e^{\epsilon}t\right) \right)dt + \int_{e^{-\epsilon}}^1 L_g \left( \tilde{\gamma}(t),\dot{\tilde{\gamma}}(t)\right) dt \\
& =& e^{\epsilon} V_{g,\Delta}(x(\epsilon)) +  \int_{e^{-\epsilon}}^1 \frac{1}{t^2} \sum_{i=1}^m \tilde{u}_i(x(-\ln t))^2 dt \\
& = &  e^{\epsilon} V_{g,\Delta}(x(\epsilon)) + \int_0^{\epsilon} e^s \sum_{i=1}^m \tilde{u}_i(x(s))^2 ds.
\end{eqnarray*}
Using (\ref{couille2}), (\ref{midi0}), and (\ref{midi}), one has, for almost every $s\in [0,\epsilon]$,
$$
 \sum_{i=1}^m \tilde{u}_i(x(s))^2 = \sum_{i=1}^m \frac{1}{4} \Big( \zeta_s \left(f_i(x(s))\right) \Big)^2 = V_{g,\Delta}(x(s)) = V_{g,\Delta}(x) e^{-2s},
 $$
and, since $V_{g,\Delta}(x(\epsilon) =  V_{g,\Delta}(x) e^{-2\epsilon}$, it follows that 
\begin{eqnarray*}
J_g(\tilde{\gamma}(\cdot)) = V_{g,\Delta}(x).
\end{eqnarray*}
Hence, $\tilde{\gamma}$ is a minimizing curve of the Bolza problem $(\mbox{BP})_{g,\Delta}$ between $\bar{x}$ and $x$. From Lemma \ref{lem6}, it cannot stay on $\mathcal{S}$ on positive times. This yields a contradiction.
\end{proof}

It follows from this lemma, and from the construction of $X$ using optimal controls, that any Carath\'eodory trajectory of $\dot{x}=X(x)$, $x(0)=x$, tends to $\bar{x}$ as $t$ tends to $+\infty$. The property of Lyapunov stability is obvious to verify. Finally, the fact
that the set $\mathcal{S}$ has Hausdorff dimension lower than or equal
to $ n-1$ is a consequence of Lemma \ref{lem4}.

%%%%%%%%%%%%%%%%%%%%%%%%%%%%%%%%%%%%%%%%%%%%%%%%%%%%%%%%%%%%%%%%%%%%%%

\subsection{Proof of Theorem \ref{THM2}}\label{secproofthm2}

Let $g$ be a Riemannian metric on $M$ and $\bar{x}$ be fixed. Since
$\Delta$ is a smooth distribution of rank two on $M$, for every $x\in
M$, there exists a neighborhood $\mathcal{V}_x$ of $x$ and two smooth
vector fields $f_1^x, f_2^x$ which represent $\Delta$ in
$\mathcal{V}_x$, that is, such that
$$
\Delta (y) = \SPAN \left\{f_1^x(y),f_2^x (y) \right\}, \quad \forall y \in \mathcal{V}_x.
$$
Moreover, as recalled in the introduction, since $\Delta$ is a
Martinet distribution, for every $x\in \Sigma_{\Delta}$, the two
vector fields $f_1^x,f_2^x$ can be chosen as
\begin{equation}\label{deff1f2}
f_1^x = \frac{\partial}{\partial x_{1}} + x_2^2\frac{\partial}{\partial
  x_{3}} \quad \mbox{and} \quad  f_2^x = \frac{\partial}{\partial x_{2}},
\end{equation}
in local coordinates. Recall that, in the neighborhood
$\mathcal{V}_x$, the Martinet surface $\Sigma_\Delta$ coincides with
the surface $x_2=0$, and the singular paths are the
integral curves of the vector field $\frac{\partial}{\partial x_{1}}$
restricted to $x_2=0$. For convenience, consider that the vector
fields $f_1^x, f_2^x$ are defined as well outside the neigborhood
$\mathcal{V}_x$. Thus, without loss of generality, for every $x\in M$,
we assume that the vector fields $f_1^x, f_2^x$ are well defined,
smooth on $M$ and satisfy
$$
 f_1^x (y) =f_2^x (y) =0, \quad \forall y \in M \setminus \mathcal{W}_x,
$$
with $\mathcal{V}_x\subset \mathcal{W}_x$, and
$$
 \SPAN \left\{f_1^x(y),f_2^x (y) \right\} \subset \Delta (y), \quad \forall y \in M.
 $$
By compactness of $\Sigma_{\Delta}$, there is a finite number of points $(x_i)_{i\in I}$ of $\Sigma_{\Delta}$ such that 
$$
\Sigma_{\Delta} \subset \cup_{i\in I} \mathcal{V}_{x_i}.
$$ 
Let $\beta :M \rightarrow [0,\infty)$ be a smooth function such that 
$$
\forall x \in M, \quad \beta (x) =0 \Longleftrightarrow x \in \Sigma_{\Delta}.
$$
For every $i\in I$, define the smooth vector field $g_i$, in local coordinates, by 
$$
g_i(y) := \beta (y) f_1^{x_i} (y), \quad \forall y \in M.
$$
By compactness of $M$, there is a finite number of points $(y_j)_{j\in J}$ of $M$ such that
$$
M \subset \left( \cup_{i\in I} \mathcal{V}_{x_i} \right)  \cup \left( \cup_{j\in J} \mathcal{V}_{y_j} \right)
$$
and 
$$
\Sigma_{\Delta} \cap \left( \cup_{j\in J} \mathcal{W}_{y_j} \right) = \emptyset.
$$
By construction, we have 
\begin{equation}
\label{quiche1}
\SPAN \left\{ g_i(y),f_2^{x_i}(y), f_1^{y_j} (y), f_2^{y_j}(y) \ \vert\ i\in I, j\in J\right\} = \Delta (y), \quad \forall y \in M \setminus \Sigma_{\Delta}
\end{equation}
and 
\begin{equation}
\label{quiche2}
\SPAN \left\{ g_i(y),f_2^{x_i}(y), f_1^{y_j} (y), f_2^{y_j}(y) \ \vert\ i\in I, j\in J\right\} \cap T_y \Sigma_{\Delta} = \{0\}, \quad \forall y \in \Sigma_{\Delta}.
\end{equation}
Indeed, for every $y \in \Sigma_{\Delta}$, there holds
\begin{equation}\label{quiche3}
\SPAN \left\{ g_i(y),f_2^{x_i}(y), f_1^{y_j} (y), f_2^{y_j}(y) \ \vert\ i\in I, j\in J\right\} = 
\SPAN \left\{ f_2^{x_i}(y)\ \vert\ i\in I\right\}.
\end{equation}
It follows from (\ref{quiche1}) and (\ref{quiche2}) that any trajectory, solution of the control system 
\begin{equation}
\label{quichecontrol}
\dot{x}(t) = \sum_{i\in I} u^1_i (t) g_i (x(t)) + u_i^2(t) f_2^{x_i}(x(t)) + \sum_{j\in J} v_j^1 (t) f_1^{y_j} (x(t)) + v_j^2(t) f_2^{y_j} (x(t)),
\end{equation} 
where $u(\cdot)=(u_1^1(\cdot),u_1^2(\cdot),\cdots, u_{\vert I\vert}^1(\cdot), u_{\vert I\vert}^2(\cdot), v_1^1(\cdot),v_1^2(\cdot), \cdots, v_{\vert J\vert}^1(\cdot), v_{\vert J\vert}^2(\cdot))$ belongs to the control set $\mathcal{U}$ defined by
$$
\mathcal{U} := L^{\infty} \left( [0,1]; \R^{2\vert I\vert+2\vert J\vert}\right),
$$
is an horizontal path of $\Delta$. Note that, for every $u(\cdot)\in \mathcal{U}$,
there exists a unique  absolutely continuous curve
$\gamma_{u(\cdot)}:[0,1]\rightarrow M$ such that
$\gamma_{u(\cdot)}(0)=\bar{x}$ and
\begin{equation*}
\begin{split}
\dot{\gamma}_{u(\cdot)}(t)
  =  & \sum_{i\in I} \left( u^1_i (t) g_i  (\gamma_{u(\cdot)}(t)) +  u_i^2(t) f_2^{x_i}(\gamma_{u(\cdot)}(t)) \right) \\
& + \sum_{j\in J} \left( v_j^1 (t) f_1^{y_j} (\gamma_{u(\cdot)}(t)) + v_j^2(t) f_2^{y_j} (\gamma_{u(\cdot)}(t))  \right) ,
\end{split}
\end{equation*}
for almost every $t\in [0,1]$. Moreover, it is clear by construction of the control system under consideration that, for every $x\in M$, there exists a control $u(\cdot) \in \mathcal{U}$ such that $\gamma_{u(\cdot)} (1) =x$. For every $u(\cdot) \in \mathcal{U}$, set 
$$
J(u(\cdot)) := \int_0^1 \left( \sum_{i\in I} \left( u_i^1(t)^2 + u_i^2(t)^2 \right) + \sum_{j\in J} \left( v_j^1(t)^2 + v_j^2(t)^2\right) \right) dt.
$$
Define the value function $W:M \rightarrow \R$ by 
$$
W(x) := \inf \left\{  J (u(\cdot)) \ \vert \ u(\cdot) \in
  \mathcal{U},\ \gamma_{u(\cdot)}(0)=\bar{x},\ \gamma_{u(\cdot)}(1)=x \right\},
$$
for every $x \in M$. By coercivity of the cost function, it is easy to
prove that, for every $x\in M\setminus \{\bar{x}\}$, there exists a
control $u(\cdot) \in \mathcal{U}$ such that
$\gamma_{u(\cdot)} (1)=x$ and
$W(x) = J(u(\cdot))$ (i.e., a minimizing control).
Moreover, by construction of the control system, more precisely, from (\ref{quiche3}), the trajectory
$\gamma_{u(\cdot)}(\cdot)$ cannot stay on the Martinet surface on a nontrivial
subinterval of $[0,1]$. As a consequence, since any singular trajectory is contained in the Martinet surface, any nontrivial minimizing
control is nonsingular. Using similar arguments as in the proof of Theorem \ref{THM1}, it follows that the value function $W$ is a viscosity solution of a certain
Hamilton-Jacobi equation, is continuous on $M$, and is locally
semiconcave in $M \setminus \{\bar{x}\}$ (see \cite{cr04}). Moreover,
the optimal trajectories of the optimal control problem under
consideration share the same properties as those of the Bolza problem
$(\mbox{BP})_{g,\Delta}$. The construction of a stabilizing
feedback then follows the same lines as in Theorem \ref{THM1}.

\begin{remark}\label{remendproof2}
For a noncompact manifold $M$, the above proof needs to be adapted by replacing a finite number of controls $(u_i)_{i\in I}$ and $(v_j)_{j\in J}$ with a \textit{locally finite} set of controls.
\end{remark}

%%%%%%%%%%%%%%%%%%%%%%%%%%%%%%%%%%%%%%%%%%%%%%%%%%%%%%%%%%%%%%%%%%%%%
%%%%%%%%%%%%%%%%%%%%%%%%%%%%%%%%%%%%%%%%%%%%%%%%%%%%%%%%%%%%%%%%%%%%%

\section{Appendix}

\subsection{Proof of Lemma \ref{LASTlem}}\label{proofofLASTlem}

Without loss of generality, we assume that $M=\R^n$. Given $k\in
\{1,\cdots ,n\}$ and $\rho >0$, denote by
$\Sigma_{\rho}^k (u)$ the set of all $x\in U$ such that $D^+u(x)$
contains a $k$-dimensional sphere of radius $\rho$, and define 
$$
\Sigma^k (u) := \left\{ x\in U \ \vert \ \mbox{dim} \left( D^+u(x)\right) =k \right\}.
$$ 
By well known properties of convex sets, one has $\Sigma^k(u) \subset \bigcup_{\rho >0}
\Sigma^k_{\rho}(u)$. Note that a point $x\in
\Sigma_{\rho}^k (u)$ does not necessarily belong to $\Sigma^k (u)$,
since $D^+ u(x)$ may be of dimension greater than $k$. The following
result is fundamental for the proof of Lemma \ref{LASTlem} (we refer
the reader to \cite{cs04} for its proof).

\begin{lemma}
\label{garance}
For every $k \in \{1,\cdots ,n\}$ and every $\rho >0$, the set
$\Sigma_{\rho}^k (u)$ is closed and satisfies
$$
\mbox{\rm Tan} \left(x,\Sigma_{\rho}^k (u) \right) \subset \left[ D^+
  u(x) \right]^{\perp}, \quad \forall x \in \Sigma^k_{\rho}(u) \cap
\Sigma^k(u)\footnote{Here, $\mbox{Tan}    \left(x,\Sigma_{\rho}^k (u)
  \right) $ denotes the tangent set to $\Sigma_{\rho}^k (u)$ at
  $x$. Recall that, given a closed set $S\subset \R^n$ and $x\in
  S$, the tangent set to $S$ at $x$, denoted by $\mbox{Tan} (x,S)$, is
  defined as the vector space generated by the set
$$
T(x,S) := \left\{ \lim_{i \rightarrow \infty} \frac{x_i -x }{t_i} \ \vert \ x_i \in S, x_i \rightarrow x, t_i \in \R_+, t_i \downarrow 0\right\}.
$$
Recall also that, if $A\subset \R^n$, then the set $A^{\perp}$ is
defined as the set of vectors $v\in \R^n$ such that $\langle
v,p\rangle =\langle v,p'\rangle $ for any $p,p' \in A$.}.
$$
\end{lemma}

Return to the proof of Lemma \ref{LASTlem}. First, note that the map $t \in
[a,b] \mapsto u(\gamma(t))$ is Lipschitzian. Hence, by
Rademacher's Theorem, it is differentiable almost everywhere on
$[a,b]$. Moreover, by the chain rule for Clarke's generalized
gradients (see \cite{clsw98}), for every $t\in [a,b]$ where $\gamma$
is differentiable, there exists $p\in \partial u(\gamma(t))$ such that
\begin{eqnarray}
\label{garance1}
\frac{d}{dt} \left( u(\gamma (t))\right) = \langle p, \dot{\gamma}(t) \rangle.
\end{eqnarray} 
For every $k\in \{1,\cdots, n\}$ and any positive integer $l$, set
$$
I_{k,l} := \left\{ t\in [a,b] \ \vert \ \gamma (t) \in \left(
    \Sigma^k_{\frac{1}{l}} (u) \cap \Sigma^k (u) \right) \setminus
  \Sigma^k_{\frac{1}{l+1}} (u) \right\}
$$
and 
$$
J := [a,b] \setminus \bigcup_{k,l} I_{k,l}.
$$
Notice that, since $u$ is locally semiconcave and $\gamma $ is
locally Lipschitzian, $u$ is differentiable at almost every $\gamma (t)$
with $t\in J$. Thus, for every such $t$, there holds necessarily
$p(t)=\nabla u(\gamma(t))$ and
\begin{eqnarray*}
\frac{d}{dt} \left( u(\gamma (t))\right) = \langle p(t), \dot{\gamma}(t) \rangle.
\end{eqnarray*} 
It remains to prove that this equality holds for almost every $t$ in
$[a,b] \setminus J$. From the Lebesgue density theorem, there exists a
sequence of measurables sets  $\{I_{k,l}^{\prime}\}$ such that all
sets $I_{k,l} \setminus I_{k,l}^{\prime}$ have Lebesgue measure zero
and such that any point in one of the sets $I_{k,l}^{\prime}$ is a
density point in that set. It is sufficient to prove the required
equality on each set $I_{k,l}^{\prime}$. Fix $k,l$ and $t \in
I_{k,l}^{\prime}$, set $x:=\gamma(t)$. Since $x$ is a density point in
$I_{k,l}^{\prime}$, there exists a sequence $\{t_i\}$ of times in
$I_{k,l}^{\prime}$ converging to $t$. Thus, the vector
$\dot{\gamma}(t)$ belongs to   $
\mbox{\rm Tan} \left(x,\Sigma_{\rho}^k (u) \right)$. Then, from Lemma
\ref{garance}, $\dot{\gamma}(t)$ belongs to  $\left[ D^+ u(x)
\right]^{\perp}$. By (\ref{garance1}), we obtain the desired
equality. This concludes the proof of Lemma \ref{LASTlem}.

\subsection{Proof of Lemma \ref{moule5}}\label{proofmoule5}
The proof that we present here is taken from \cite{riffordmemoir}
(compare with \cite{ls95,op01}). For the sake of simplicity, assume that $M=\R^n$, endowed with
the Euclidean metric. Since the property to be proved is local,
we assume that there are $m$ smooth
vector fields $f_1,\cdots, f_m$, orthonormal with respect to
the Euclidean metric, such that
$$
\Delta (x) = \SPAN \left\{ f_i(x) \ \vert \ i=1,\ldots ,m
\right\},
$$
in a neighborhood $\mathcal{V}$ of $\bar{x}$. With these notations,
the associated Hamiltonian
%% \footnote{As the reader can see, we prefer to work
%% momentarily with an Hamitonian defined on $\R^n \times \R^n$ than on
%% $\R^n \times \left(\R^n \right)$}
$H: \R^n \times \R^n \rightarrow
\R$ is
$$
H(x,p) := \max_{u\in \R^m} \left\{ \langle p, \sum_{i=1}^m
  u_if_i(x)\rangle  - \frac{1}{2} \sum_{i=1}^m u_i^2 \right\} =
\frac{1}{2} \sum_{i=1}^m \langle p,f_i(x)\rangle^2,
$$
for every $(x,p) \in R^n \times \R^n$.\\
Our aim is now to prove the following result: for every $p_0 \in \R^n$
such that $H(\bar{x},p_0)\neq 0$, there exist a neighborhood
$\mathcal{W}$ of $p_0$ in $\R^n $ and $\epsilon >0$ such that every
solution $(x(\cdot),p(\cdot)):[0,\epsilon] \rightarrow \R^n \times
\R^n$ of the Hamiltonian system
\begin{eqnarray}
 \label{swen1}
 \left\{
 \begin{array}{lll}
\displaystyle
 \dot{x}(t) & = & \displaystyle\frac{\partial H}{\partial p} (x(t),p(t)) = \displaystyle\sum_{i=1}^m \langle p(t),f_i(x(t))\rangle f_i(x(t)) \\
\displaystyle \dot{p}(t) & = & - \displaystyle\frac{\partial H}{\partial x} (x(t),p(t)) =\displaystyle
 -\sum_{i=1}^m \langle p(t),f_i(x(t))\rangle df_i(x(t))^* p(t),
 \end{array}
 \right.
 \end{eqnarray}
with $x(0)=\bar{x}$ and $p(0)\in \mathcal{W}$, satisfies
\begin{eqnarray}
\label{swen2}
\int_0^{\epsilon} \sum_{i=1}^m \langle p(t),f_i(x(t))\rangle^2 dt \leq \int_0^{\epsilon} \sum_{i=1}^m u_i(t)^2 dt ,
\end{eqnarray}
for every control $u(\cdot) \in L^{\infty}([0,\epsilon];\R^m)$ such that the solution of 
\begin{eqnarray}
\label{swen3}
\dot{y}(t) = \sum_{i=1}^m u_i(t) f_i(y(t)), \quad y(0)=\bar{x},
 \end{eqnarray}
satisfies $y(\epsilon ) =x(\epsilon)$. Let $p_0 \in \R^n \setminus \{0\}$ be fixed,  we need the following lemma.

\begin{lemma}
\label{lemswen}
There exist a neighborhood $\mathcal{W}$ of $p_0$ and $\rho >0$ such
that, for every $p\in \mathcal{W}$, there exists a function $S:
B(\bar{x},\rho) \rightarrow \R$ of class $C^1$ which satisfies
\begin{eqnarray}
\label{swen4}
H(x,\nabla S(x)) = H(\bar{x},p), \quad \forall x \in B(\bar{x},\rho),
\end{eqnarray}
and such that, $(x^p(\cdot),p^p(\cdot))$ denotes the solution of
(\ref{swen1}) satisfying $x^p(0)=\bar{x}$ and $p^p(0)=p$, then
\begin{eqnarray}
\label{swen5}
\nabla S(x^p(t))=p^p(t), \quad \forall t\in (-\rho,\rho).
\end{eqnarray}
\end{lemma}

\begin{proof}
The proof consists in applying the method of characteristics. Let
$\Pi$ be the linear hyperplane such that $\langle p_0,v\rangle =0$ for
every $v\in \Pi$. We first show how to construct locally
$S$ as the solution of the Hamilton-Jacobi equation (\ref{swen4})
which vanishes on $\bar{x} + \Pi$ and such that $\nabla
S(\bar{x})=p_0$.
Up to considering a smaller neighborhood $\mathcal{V}$,
we assume that $H(x,p_0)\neq 0$ for every $x\in
\mathcal{V}'$. For every $x\in (\bar{x}+\Pi) \cup \mathcal{V}$, set
$$
\bar{p}(x) := \sqrt{\frac{H(\bar{x},p_0)}{H(x,p_0)}} p_0.
$$
Then, $H(x,\bar{p}(x))=H(\bar{x},p_0)$ and $\bar{p}(x) \perp \Pi$, for
every $x\in \mathcal{V}'$. There exists $\mu >0$ such that, for every
$x\in (\bar{x}+\Pi) \cup \mathcal{V}$, the solution
$(x_x(\cdot),p_x(\cdot))$ of (\ref{swen1}), satisfying $x_x(0)=x$ and
$p_x(0)=\bar{x}$, is defined on the interval $(-\mu,\mu)$. For
every $x\in  (\bar{x} +\Pi) \cup \mathcal{V}$ and every $t\in
(-\mu,\mu)$, set $\theta(t,x) := x_x(t)$. The mapping $(t,x) \mapsto
\theta(t,x)$ is smooth. Moreover, $\theta(0,x)=x$ for every $x\in
(\bar{x} +\Pi) \cup \mathcal{V}$ and $\dot{\theta}(0,\bar{x})=
\sum_{i=1}^m \langle \bar{p}(x),f_i(\bar{x})\rangle f_i(\bar{x}) $
does not belong to $\Pi$. Hence there exists $\rho \in (0,\mu)$ with
$B(\bar{x},\rho) \subset \mathcal{V}$ such that the mapping $\theta$
is a smooth diffeomorphism from $(-\rho,\rho) \times \left(
  (\bar{x}+\Pi) \cup B(\bar{x},\rho) \right)$ into a neighborhood
$\mathcal{V}'$ of $\bar{x}$. Denote by $\varphi = (\tau,\pi)$ the
inverse function of $\theta$, that is the function such that $(\theta \circ \varphi)(x)=(\tau(x),\pi(x))=x$ for every $x\in \mathcal{V}'$. Define the two vector fields $X$ and $P$ by 
$$
X(x) := \dot{\theta} (\tau(x),\pi(x)) \quad \mbox{and} \quad P(x) :=  p_{\pi(x)}(\tau(x)), \quad \forall x\in \mathcal{V}'.
$$
Then,
\begin{eqnarray*}
X(\theta(t,x))= \dot{\theta}(t,x) = \dot{x}_x(t) & = & \sum_{i=1}^m \langle p_x(t),f_i(x_x(t)) \rangle f_i(x_x(t)) \\
& = &\sum_{i=1}^m \langle P(\theta(t,x)), f_i(\theta(t,x))\rangle f_i(\theta(t,x)),
\end{eqnarray*}
and 
$$
\sum_{i=1}^m \langle P(\theta(t,x)), f_i(x_x(t))\rangle^2 = \sum_{i=1}^m \langle p_x(t),f_i(x_x(t))\rangle^2 = 2 H(x,\bar{p}(x))   = 2 H(\bar{x},p_0),
$$
for every $t \in (-\rho,\rho)$ and every $x \in  (\bar{x}+\Pi) \cup
B(\bar{x},\rho)$. For every $x\in \mathcal{V}'$, set $\alpha_i(x) :=
\langle P(x),f_i(x)\rangle$. Hence,
$$
X(x)= \sum_{i=1}^m \alpha_i(x) f_i(x) \quad \mbox{and}  \quad
\sum_{i=1}^m \alpha_i(x)^2 = H(\bar{x},p_0),
$$
for every $x\in \mathcal{V}'$. Define the function $S:\mathcal{V}'
\mapsto \R$ by
$$
S(x):=2H(\bar{x},p_0)\tau(x), \quad \forall x\in \mathcal{V}'.
$$ 
We next prove that $\nabla S(x) = P(x)$ for every $x\in \mathcal{V}'$.
For every $t\in (-\rho,\rho)$, denote by
$W_t:=\{y \in \mathcal{V}'\ \vert\ \tau(y)=t\}$.
In fact, $W_t$ coincides with the
set of $y\in \mathcal{V}'$ such that $S(y)=2H(\bar{x},p_0) t$. It is a
smooth hypersurface which satisfies $\nabla S(y) \perp T_yW_t$ for
every $y \in W_t$. Let $y\in W_t$ be fixed, there exists $x\in
(\bar{x}+\Pi) \cup B(\bar{x},\rho)$ such that
$y=\theta(t,x)=x_x(t)$. Let us first prove that $P(y)=p_x(t)$ is
orthogonal to $T_yW_t$. To this aim, without loss of
generality we assume that $t>0$. Let $w\in T_yW_t$, there exists $v\in
\Pi$ such that $w=d_x\theta_t(x)v$. For every $s \in [0,t]$, set
$z(s):=d_x\theta(s,x)v$. We have
\begin{eqnarray*}
\dot{z}(s) =\frac{d}{ds} d_x\theta (s,x)v  = \frac{d}{dx} \dot{\theta}(t,x) v = \frac{d}{dx} X(\theta(t,x)) v  =  dX(\theta(t,x)) z (s).
\end{eqnarray*}
Hence,
\begin{eqnarray*}
\frac{d}{ds} \langle z(s), p_x(s) \rangle & = & \langle \dot{z}_(s),p_x(s) \rangle + \langle z(s), \dot{p}_x(s) \rangle \\
& = & \langle dX(\theta(s,x)) z(s), p_x(s) \rangle - \langle z(s), \sum_{i=1}^m \langle p_x(s),f_i(x_x(s)) \rangle df_i(x_x(s))^* p_x(s) \rangle.
\end{eqnarray*}
Since $X(x) = \sum_{i=1}^m \alpha_i(x) f_i(x)$ and $\sum_{i=1}^m
\alpha_i(x)^2 = H(\bar{x},p_0)$ for every $x\in \mathcal{V}'$, there
holds
\begin{eqnarray*}
dX(x_x(s))^*p_x(s) & = & \sum_{i=1}^m \alpha_i(x_x(s)) df_i(x_x(s))^*p_x(s)   + \sum_{i=1}^m \langle f_i(x_x(s)),p_x(s)\rangle \nabla \alpha_i (x_x(s))\\
& = & \sum_{i=1}^m \alpha_i(x_x(s)) df_i(x_x(s))^*p_x(s)   + \sum_{i=1}^m \alpha_i(x_x(s)) \nabla \alpha_i(x_x(s)) \\
& = &  \sum_{i=1}^m \alpha_i(x_x(s)) df_i(x_x(s))^*p_x(s).
\end{eqnarray*}
We deduce that $\frac{d}{ds} \langle z(s), p_x(s) \rangle =0$ for every $s\in [0,t]$. Hence,
$$
\langle w, P(y)\rangle = \langle w,p_x(t) \rangle=\langle z(t),p_x(t)\rangle = \langle z(0), \bar{p}(x)\rangle = 0.
$$
This proves that $P(y)$ is orthogonal to $T_yW_t$, which implies that
$P(y)$ and $\nabla S(y)$ are collinear. Furthermore, since
$S(x_x(s)) =2H(\bar{x},p_0) s$ for every $s\in [0,t]$, one gets
$$
\langle \nabla S(x_x(t)), \dot{x}_x(t) \rangle = 2H(\bar{x},p_0) =  \langle p_x(t), \dot{x}_x(t) \rangle. 
$$
Since $\dot{x}_x(t)=X(y)$ does not belong to $T_yW_t$, we deduce that $\nabla S(x_x(t)) = p_x(t)$. In consequence, we proved that $\nabla S(x)= P(x)$ for every $x\in \mathcal{V}'$.
\end{proof}

Let us now conclude the proof of Lemma \ref{moule5}. Clearly, there
exists $\epsilon >0$ such that
every solution  $(x(\cdot),p(\cdot)):[0,\epsilon] \rightarrow \R^n
\times \R^n$ of (\ref{swen1}), with $x(0)=\bar{x}$ and $p(0)\in
\mathcal{W}$, satisfies
$$
x(t) \in B(\bar{x},\rho), \quad \forall t \in [0,\epsilon].
$$
Moreover, we have 
$$
 S(x(\epsilon)) - S(\bar{x}) = 2\epsilon H(\bar{x},p).
$$
Let $u(\cdot) \in L^{\infty}([0,\epsilon];\R^m)$ be a control such
that the solution $y(\cdot):[0,\epsilon] \rightarrow \mathcal{W}$ of
(\ref{swen3}) starting at $\bar{x}$ satisfies $y(\epsilon )
=x(\epsilon)$. We have
\begin{eqnarray*}
S(x(\epsilon)) - S(\bar{x}) & =  & S(y(\epsilon)) - S(y(0)) \\
& = & \int_0^{\epsilon} \frac{d}{dt} \left( S(y(t)) \right) dt \\
& = &\int_0^{\epsilon} \langle \nabla S(y(t)), \dot{y}(t) \rangle dt \\
& \leq & \int_0^{\epsilon} H(y(t),dS(y(t))) + \frac{1}{2} \sum_{i=1}^m u_i(t)^2 dt \\
& = & \epsilon H(\bar{x},p) + \int_0^{\epsilon} \sum_{i=1}^m u_i(t)^2 dt.
\end{eqnarray*}
The conclusion follows.

%%%%%%%%%%%%%%%%%%%%%%%%%%%%%%%%%%%%%%%%%%%%%%%%%%%%%%%%%%%%%%%%%%%%%%%%%%%%%
%%%%%%%%%%%%%%%%%%%%%%%%%%%%%%%%%%%%%%%%%%%%%%%%%%%%%%%%%%%%%%%%%%%%%%%%%%%%%


\begin{thebibliography}{99}

\bibitem{agrachev98}
A.~Agrachev.
\newblock Compactness for Sub-Riemannian length-minimizers and
subanalyticity.
\newblock Control theory and its applications (Grado, 1998), {\em
  Rend.\ Sem.\ Mat.\ Univ.\ Politec.\ Torino}, 56(4):1--12, 2001.

\bibitem{ABCK}
A.~Agrachev, B.~Bonnard, M.~Chyba, and I.~Kupka.
\newblock Sub-Riemannian sphere in Martinet flat case.
\newblock {\em ESAIM Control Optim.\ Calc.\ Var.} 2:377--448, 1997.

\bibitem{asach}
A.~Agrachev and Y.~Sachkov.
\newblock {\em Control theory from the geometric viewpoint}.
\newblock Encyclopaedia of Mathematical Sciences, 87,
Control Theory and Optimization, II, Springer-Verlag, Berlin, 2004.

\bibitem{AS}
A.~Agrachev and A.~Sarychev.
\newblock Sub-Riemannian metrics: minimality of
singular geodesics versus subanalyticity.
\newblock {\em ESAIM Control Optim. Calc. Var.}, 4:377--403, 1999.

%\bibitem{as04}
 %A.~A.~Agrachev and Y.~L.~Sachkov.
%\newblock {\em Control theory from the geometric viewpoint}.
%\newblock Encyclopaedia of Mathematical Sciences, 87. Springer-Verlag, Berlin, 2004.

\bibitem{aac92}
G.~Alberti, L.~Ambrosio and P.~Cannarsa.
\newblock On the singularities of convex functions.
\newblock {\em Manuscripta Math.}, 76(3-4):421--435, 1992.

\bibitem{anbr99}
F.~Ancona and A.~Bressan.
\newblock Patchy vector fields and asymptotic stabilization.
\newblock {\em ESAIM Control Optim. Calc. Var.}, 4:445--471, 1999.

\bibitem{barles94}
G.~Barles.
\newblock {\em Solutions de viscosit\'e des \'equations de Hamilton-Jacobi}.
\newblock Math\'ematiques \& Applications, 17. Springer-Verlag, Berlin, 1994.

\bibitem{bcd97}
M.~Bardi and I.~Capuzzo-Dolcetta.
\newblock {\em Optimal control and viscosity solutions of
  {H}amilton-{J}acobi-{B}ellman equations}.
\newblock Systems \& Control: Foundations \& Applications. Birkh\"auser Boston Inc., Boston, MA, 1997.

\bibitem{bellaiche96}
A.~Bella\"iche.
\newblock The tangent space in sub-Riemannian geometry.
\newblock in {\em Sub-Riemannian Geometry, Birkh\"auser}, 1--78, 1996. 

\bibitem{Bismut}
J.-M.~Bismut.
\newblock {\em Large deviations and the Malliavin
calculus}.
\newblock Progress in Mathematics 45, Birkh\"{a}user, 1984.

%\bibitem{BC}
%B.~Bonnard and M.~Chyba.
%\newblock {\em Singular trajectories and their role in control theory}.
%\newblock Math.\ \& Appl.\ 40, Springer-Verlag, 2003.

\bibitem{BTtoulouse}
B.~Bonnard and E.~Tr\'elat.
\newblock On the role of abnormal minimizers in sub-Riemannian
geometry.
\newblock {\em Ann. Fac. Sci. Toulouse Math. (6)}, 10(3):405--491, 2001.

\bibitem{bredon93}
G.~E.~Bredon.
\newblock {\em Topology and Geometry}.
\newblock Graduate Texts in Mathematics, vol. 139. Springer-Verlag, New York,  1993.

%\bibitem{Bressantorino}
% A.~Bressan.
%\newblock Singularities of stabilizing feedbacks.
%\newblock {\em Rend. Sem. Mat. Univ. Politec. Torino}, 56(4):87--104, 1998.

\bibitem{brockett}
R.~W.~Brockett.
\newblock {\em Asymptotic stability and feedback stabilization}.
\newblock Differential geometric control theory, R.~W.~Brockett,
R.~S.~Millman and H.~J.~Sussmann, ed., Boston, Birkh\"auser,
181--191, 1983.

%\bibitem{BryantHsu}
%R.~L.~Bryant and L.~Hsu.
%\newblock Rigidity of integral curves of rank 2 distributions.
%\newblock {\em Invent.\ Math.} 114:435--461,1993.

\bibitem{cr04}
P.~Cannarsa and L.~Rifford.
\newblock Semiconcavity results for optimal control problems admitting
no singular minimizing controls.
\newblock Preprint, 2006.
    
\bibitem{cs04}
P.~Cannarsa and C.~Sinestrari.
\newblock {\em Semiconcave functions, Hamilton-Jacobi equations, and
optimal control}.
\newblock Progress in Nonlinear Differential Equations and their
Applications, 58. Birkh\"auser Boston Inc., Boston, MA, 2004.

\bibitem{CJTcras}
Y.~Chitour, F.~Jean and E.~Tr\'elat.
\newblock Propri\'et\'es g\'en\'eriques des trajectoires
singuli\`eres.
\newblock {\em C. R. Acad. Sci. Paris S\'er. I Math.}, 337(1):49--52, 2003.

\bibitem{CJT}
Y.~Chitour, F.~Jean and E.~Tr\'elat.
\newblock Genericity results for singular curves.
\newblock  {\em J.\ Diff.\ Geom.} 73(1):45--73, 2006.

\bibitem{CJTnew}
Y.~Chitour, F.~Jean and E.~Tr\'elat.
\newblock Singular trajectories of control-affine systems.
\newblock Accepted for publication in SIAM J. Cont. Optim. 

\bibitem{chow39}
W.~L.~Chow.
\newblock \"Uber Systeme von linearen partiellen
Differentialgleichungen ester Ordnung.
\newblock {\em Math. Ann.}, 117: 98--105, 1939.

\bibitem{clarke}
F.~H.~Clarke.
\newblock Optimization and nonsmooth analysis.
\newblock Second Edition. {\em Classics in Applied Math.} 5, SIAM,
Philadelphia, 1990.

\bibitem{clss97}
F.~H. Clarke, Yu.~S. Ledyaev, E.D. Sontag, and A.I. Subbotin.
\newblock Asymptotic controllability implies feedback stabilization.
\newblock {\em I.E.E.E. Trans. Aut. Control}, 42:1394--1407, 1997.

\bibitem{clsw98}
F.~H.~Clarke, Yu.~S.~Ledyaev, R.~J.~Stern and P.~R.~Wolenski.
\newblock {\em Nonsmooth Analysis and Control Theory}.
\newblock Graduate Texts in Mathematics, vol. 178. Springer-Verlag, New York,   1998.

\bibitem{coron90}
J.-M.~Coron.
\newblock A necessary condition for feedback stabilization. 
\newblock {\em Systems Control Lett.}, 14(3):227--232, 1994.

\bibitem{coron92}
J.-M.~Coron.
\newblock Global asymptotic stabilization for controllable systems
without drift.
\newblock {\em Math. Control Signals Systems}, 5(3):295--312, 1992.

%\bibitem{coron94}
% J.-M.~Coron.
%\newblock Linearized control systems and application to smooth stabilization. 
%\newblock {\em SIAM J. Control Optim.}, 32(2):358--386, 1994.

\bibitem{CL} 
M.~G.~Crandall and P.-L.~Lions.
\newblock Viscosity solutions of Hamilton-Jacobi equations.
\newblock {\em Trans.\ Amer.\ Math.\ Soc.}, 277(1):1--42, 1983.

\bibitem{federer69}
H.~Federer.
\newblock {\em Geometric measure theory}.
\newblock Die Grundlehren der mathematischen Wissenschaften, Band 153.
Springer-Verlag, New York, 1969.

\bibitem{Hardt}
R.~M.~Hardt.
\newblock Stratification of real analytic mappings and images.
\newblock {\em Invent.\ Math.}, 28:193--208, 1975.

\bibitem{Hironaka} 
H.~Hironaka.
\newblock Subanalytic sets.
\newblock In {\em Number theory , algebraic geometry and commutative algebra, in honor of Y. Akizuki.} Kinokuniya, Tokyo, 1973.

\bibitem{Hsu}
L.~Hsu.
\newblock Calculus of variations via the Griffiths formalism.
\newblock {\em J.\ Diff.\ Geom.} 36, 1992.


%\bibitem{kupka97}
% I.~Kupka.
%\newblock G\'eom\'etrie sous-riemannienne.
%\newblock {\em Ast\'erisque}, (241):Exp.\ No.\ 817, 5, 351--380, 1997.
%\newblock S\'eminaire Bourbaki, Vol.\ 1995/96.

\bibitem{ln05}
Y.~Y.~Li and L. Nirenberg.
\newblock The distance function to the boundary, Finsler geometry and the singular set   of viscosity solutions of some Hamilton-Jacobi equations.
\newblock {\em Comm. Pure Appl. Math.}, 58:85--146, 2005.

\bibitem{Lions} 
P.-L.~Lions.
\newblock {\em Generalized solutions of Hamilton-Jacobi
equations.}
\newblock Pitman (Advanced Publishing Program), Boston, Mass., 1982.

\bibitem{ls95}
W.~Liu and H.~J.~Sussmann.
\newblock Shortest paths for sub-Riemannian metrics on rank-2 distributions.
\newblock {\em Mem. Amer. Math. Soc.} 118(564), 1995. 

\bibitem{montgomery02}
R.~Montgomery.
\newblock {\em A tour of subriemannian geometries,
their geodesics and applications}.
\newblock {\em Mathematical Surveys and Monographs}, Vol.\ 91.
\newblock American Mathematical Society, Providence, RI, 2002.

\bibitem{op01}
P.~Orro and F.~Pelletier.
\newblock Differential properties of the distance function associated to a submanifold in sub-Riemannian geometry.
\newblock Pr\'epublication du LAMA, 2001.

\bibitem{P}
L.~Pontryagin, V.~Boltyanskii, R.~Gamkrelidze, and E.~Mischenko.
\newblock {\em The mathematical theory of optimal
processes}.
\newblock Wiley Interscience, 1962.

\bibitem{rashevsky38}
P.~K.~Rashevsky.
\newblock About connecting two points of a completely nonholonomic space by admissible curve.
\newblock {\em Uch.\ Zapiski Ped.\ Inst.\ Libknechta}, 2:83--94, 1938.

\bibitem{rifford02}
L.~Rifford.
\newblock Semiconcave control-Lyapunov functions and stabilizing feedbacks.
\newblock  {\em SIAM J. Control Optim.}, 41(3):659--681, 2002.

\bibitem{riffordIHP}
L.~Rifford.
\newblock Stratified semiconcave control-Lyapunov functions and the stabilization problem.
\newblock {\em Ann. Inst. H. Poincar\'e Non Lin\'eaire}, 22(3):343--384, 2005.


\bibitem{riffordMCT03}
L.~Rifford.
\newblock The Stabilization Problem: AGAS and SRS Feedbacks.
\newblock In {\em Optimal Control, Stabilization, and Nonsmooth Analysis}, 
Lectures Notes in Control and Information Sciences, 301,
Springer-Verlag, Heidelberg (2004), 173--184.

\bibitem{riffordTORINO}
L.~Rifford.
\newblock The stabilization problem on surfaces.
\newblock {\em Rend. Semin. Mat. Torino} 64(1):55--61, 2006.

\bibitem{riffordSRS}
L.~Rifford.
\newblock On the existence of local smooth repulsive stabilizing
feedbacks in dimension three.
\newblock {\em J. Differential Equations}, 226(2):429--500, 2006.

\bibitem{riffordmemoir}
L.~Rifford.
\newblock {\em Nonholonomic Variations}.
\newblock Monograph, in progress.

\bibitem{rt05}
L.~Rifford and E.~Tr\'elat.
\newblock Morse-Sard type results in sub-Riemannian geometry.
\newblock {\em Math.\ Ann.}, 332(1):145--159, 2005.

\bibitem{ryan94}
E.~P.~Ryan.
\newblock On Brockett's condition for smooth stabilizability ans its
necessity in a context of nonsmooth feedback.
\newblock {\em SIAM J. Control Optim.}, 32:1597--1604, 1994.

\bibitem{Sar} 
A.~Sarychev.
\newblock The index of the second variation of a control system.
\newblock {\em Math.\ USSR Sbornik},  41(3):383--401,1982.

\bibitem{sussmann79}
H.~J.~Sussmann.
\newblock Subanalytic sets and feedback control.
\newblock {\em J. Differential Equations}, 31(1):31--52, 1979.

\bibitem{trelat00}
E.~Tr\'elat.
\newblock Some properties of the value function and its level sets for affine control systems with quadratic cost.
\newblock {\em J.\ Dynam.\ Cont.\ Syst.}, 6(4):511--541, 2000.

\bibitem{Zh}
M.~Zhitomirski.
\newblock {\em Typical singularities of differential 1-forms and pfaffian
equations}.
\newblock Trans.\ Math.\ Monographs, Vol.\ 113, American Mathematical Society, 1992.

\end{thebibliography}
\end{document}